\input amstex
\documentstyle{amsppt}
\magnification 1100
\input xy
\xyoption{all} 
\pagewidth{420pt}
\define\a{\alpha}\define\be{\beta}\define\Gam{\varGamma}
\define\Ph{\varphi}\define\z{\zeta}\define\zn{\zeta_{p^n}}
\define\G{\text{\rm Gal}}\define\Dc{\bold D_{\text{\rm cris}}}\define\La{\varLambda}
\define\Dd{\bold D_{\text{\rm dR}}}

\define\Drig{\bold D_{\text{\rm rig}}}
\define\RG{\bold R\Gamma}
\define\res{\text{\rm res}}
\define \R{\bold R}
\define \boB{\bold B}
\define \Brigplus{\widetilde {\bold B}_{\text{\rm rig}}^+}
\define \A{\bold A}
\define \E{\bold E}
\define \bD{\bold D}
\define \co{\text{\rm c}}
\define \cl{\text{\rm cl}}
\define \Tr{\text{\rm Tr}}
\define \Tot{\text{\rm Tot}}
\define \Iw{\text{\rm Iw}}
\define \Zp{\Bbb Z_p}

\define\ep{\varepsilon}
\define\Bd{\bold B_{\text{\rm dR}}}
\define\Bc{\bold B_{\text{\rm cris}}}\define\Exp{\text{\rm Exp}}

\define\Hom{\text {\rm Hom}}
\define\Ext{\text {\rm Ext}}\define\M {\Cal M}

\define\F{\text{\rm Fil}}
\define \g{\gamma}
\define \gn{\gamma_n}

\define\Ind{\text{\rm Ind}}

\define\Hi{H_{\text{\rm Iw}}}

\define\Ddagrig{\bold D^{\dagger}_{\text{\rm rig}}}

\define\Brigdag{\CR (K)}

\define\iso{\overset{\sim}\to{\rightarrow}}

\define\Rep{\text{\bf Rep}}
\define\Gal{\text{\rm Gal}}
\define\dop{\partial}
\define \HG{\Cal H(\Gamma)}
\define \bExp{\text{\bf E}\text{\rm xp}}
\define \Vc{\bold V_{\text{\rm cris}}}
\define \CR{\Cal R}
\define \sha{\hbox{\text{\bf III}\hskip -11pt\vrule width10pt
depth0pt height0.4pt\hskip 1pt}}
\topmatter \nologo \TagsOnRight
\NoBlackBoxes
\title
On trivial zeros of Perrin-Riou's $L$-functions 
\endtitle
\author
{ Denis Benois}
\endauthor
\date April 2009
\enddate
\subjclass \nofrills 2000 Mathematics Subject Classification. 11F80, 11R23,  11G40, 11S  
\endsubjclass
\abstract In the previous paper \cite{Ben2} we generalized Greenberg's construction
of the $\Cal L$-invariant to semistable $p$-adic representations. Here we prove that
this construction is compatible with Perrin-Riou's theory of $p$-adic $L$-functions.
Namely, using   Nekov\'a\v r's  machinery of Selmer complexes we prove that
our $\Cal L$-invariant appears as an additional factor in the Bloch-Kato type formula 
for special values of Perrin-Riou's Iwasawa $L$-function.
\endabstract
\address Denis Benois, Institut de Math\'ematiques,
Universit\'e Bordeaux I, 351, cours de la Lib\'eration 33405
Talence, France
\endaddress
\email denis.benois\@math.u-bordeaux1.fr
\endemail
\leftheadtext{Denis Benois} \rightheadtext{Trivial zeros} \toc
\nofrills {\bf Table of contents} \widestnumber \head{\S 5.}
\widestnumber\subhead{5.2} \specialhead{}Introduction
\endspecialhead 
\head \S 1. $(\Ph,\Gamma)$-modules
\endhead
\subhead 1.1. Rings of $p$-adic periods 
\endsubhead
\subhead 1.2. $(\Ph,\Gamma)$-modules
\endsubhead
\subhead 1.3. Cohomology of $(\Ph,\Gamma)$-modules
\endsubhead
\subhead 1.4. Crystalline representations
\endsubhead
\head \S 2. The exponential map 
\endhead
\subhead 2.1. The Bloch-Kato exponential map
\endsubhead
\subhead 2.2. The large exponential map
\endsubhead
\head \S 3. The $\Cal L$-invariant
\endhead
\subhead 3.1. Definition of the $\Cal L$-invariant
\endsubhead
\subhead 3.2. The Bockstein homomorphism
\endsubhead
\head \S 4. Special values of $p$-adic $L$-functions 
\endhead
\subhead 4.1. The Bloch-Kato conjecture
\endsubhead
\subhead 4.2. The complex $\RG_{\text{\rm Iw},h}^{(\eta_0)}(D,V)$
\endsubhead
\subhead 4.3. The module of $p$-adic $L$-functions 
\endsubhead
\head  {Appendix. Galois cohomology of $p$-adic representations}
\endhead
\endtoc
\endtopmatter
\document

\head{\bf Introduction}
\endhead

\flushpar
{\bf 0.1.} In \cite{Ben2}, using ideas of  Colmez \cite{C4}  
we defined   a natural generalization
of Greenberg's $\Cal L$-invariant \cite{G} to pseudo-geometric representations $V$
of $\text{\rm Gal} (\overline {\Bbb Q}/\Bbb Q)$ which are semistable at $p.$
More precisely, assume that $V$ satisfies the following conditions:

1)  $H^0(V)=H^0(V^*(1))=0$ and $H^1_f(V)=H^1_f(V^*(1))=0$;

2) $V$ is semistable  at $p$ and the  map $1-p^{-1}\Ph^{-1}$ acts semisimply on $\bD_{\text{\rm st}}
(V)$.

3)  $\bD_{\text{\rm st}} (V)^{\Ph=1}=0.$

4) The $(\Ph,\Gamma)$-module $\Ddagrig (V)$ has no
crystalline subquotient of the form
$$
0@>>>\Cal R (\vert x\vert x^k)@>>>U@>>>\Cal R@>>>0, \qquad
k\geqslant 1.
$$
See sections 1.1, 2.1 and 3.1  for unexplained notations and further details. 
Remark that  $\Cal R$ denotes the Robba ring over $\Bbb Q_p$ and 4) is a 
direct generalization of Hypothesis U of \cite{G}.
Let $t_V(\Bbb Q_p)=\bD_{\text{\rm st}}(V)/\F^0 \bD_{\text{\rm st}}(V)$ denote
the tangent space of $V$ at $p$. We say that a $\Bbb Q_p$-subspace $D\subset
\bD_{\text{\rm st}} (V)$ is admissible if it is stable under the action of $\Ph$
and the natural projection $D @>>>t_V(\Bbb Q_p)$ is an isomorphism. The main construction
of \cite{Ben2} associates to $(V,D)$ a $p$-adic number $\Cal L (V,D)\in \Bbb Q_p$
which coincides with the Greenberg's $\Cal L$-invariant if $V$ is ordinary at $p$ and
$D=\bD_{\text{\rm st}} (F^1V)$ where $F^1V$ denotes the  canonical filtration of $V$
provided by ordinarity. 
\newline
\,

{\bf 0.2.} The goal of the present paper is to  show that this definition is compatible with 
Perrin-Riou's theory of $p$-adic $L$-functions. 
For a profinite group $G$ and a continuous $G$-module $X$ we denote by $C_c^\bullet (G,X)$
the standard complex of continuous cochains. Let $S$ be a finite set of primes
containing $p$. Denote by  $G_S$ the Galois group of the maximal algebraic extension of $\Bbb Q$
unramified outside $S\cup \{\infty\}$.  Set
$\RG_S(X)= C^{\bullet}_c (G_S,X)$ and $\RG (\Bbb Q_v,X)=C_c^\bullet (G_v,X)$, where 
 $G_v$ is the absolute Galois group  of  $\Bbb Q_v$. Let $\RG_{c}(V)$ denote the complex
 sitting in the distinguished triangle
$$
\RG_{c}(V) @>>>\RG_S(V)@>>>\underset{v\in S\cup\{\infty\}}\to \oplus \RG (\Bbb Q_v,V).
$$ 
The Euler-Poincar\'e line of $V$ is defined by $\Delta_{\text{\rm EP}}(V)={\det}_{\Bbb Q_p}^{-1}\RG_{c}(V).$
\newline
Now assume that $V$ is the $p$-adic realization of a pure motive $M/\Bbb Q$. 
Let $M_B$ and $M_{\text{\rm dR}}$ denote the Betti and the de Rham realizations of
$M$ and let $t_M(\Bbb Q)=M_{\text{\rm dR}}/\F^0M_{\text{\rm dR}}$ denote the tangent space of $M.$
Fixing  non zero elements $\omega_B\in {\det}_{\Bbb Q}M_B^+$ and $\omega_t\in {\det}_{\Bbb Q}t_M(\Bbb Q)$
one can define a canonical  trivialization 
$$
i_{\omega_t,\omega_B,p}\,\,:\,\,\Delta_{\text{\rm EP}}(V) @>>>\Bbb Q_p \,.
$$ 
Let $T$ be a $G_S$-stable lattice of $V.$ According to the conjecture of Bloch and Kato  \cite{BK}
in the form of Fontaine and Perrin-Riou \cite{F3}
$$
i_{\omega_t,\omega_B,p} (\Delta_{\text{\rm EP}}(T))=\frac{L(M,0)}{\Omega_\infty (\omega_t,\omega_B)}\Bbb Z_p,
$$
where $\Omega_\infty (\omega_t,\omega_B)$ is the Deligne period.
Assume in addition that  $V$  is crystalline at $p$. Fix an admissible subspace $D$ of $\Dc (V)$ and
a $\Bbb Z_p$-lattice $N$ of $D.$  From the semisimplicity of $\Ph$ we deduce the decomposition
$D\simeq D_{-1}\oplus D^{\Ph=p^{-1}}$ where $D_{-1}=(\Ph-p^{-1})\,D.$
Set $\Gamma =\text{\rm Gal}(\Bbb Q(\zeta_{p^\infty})/\Bbb Q),$
$\Gamma_1=\text{\rm Gal}(\Bbb Q(\zeta_{p^\infty})/\Bbb Q(\zeta_p))$ and $\Lambda =\Bbb Z_p[[\Gamma_1]].$
Fix a topological generator $\gamma_1\in \Gamma_1$ and denote by $\Cal H$ the ring
of operators $f(\gamma_1-1)$ where $f(X)=\sum_{n=0}^\infty a_nX^n\in \Bbb Q_p[[X]]$ converges on the $p$-adic open unit disk.
Let $\Cal K$ be the field of fractions of $\Cal H.$ Fix $h\geqslant 1$ such that $\F^{-h}\Dc (V)=\Dc (V).$
Perrin-Riou's theory \cite{PR2} associates to
$(T,N)$ a free $\Lambda$-module $\bold L^{(\eta_0)}_{\text{\rm Iw},h}(N,T)\subset \Cal K$
Fix a generator $f(\gamma_1-1) $ of  $\bold L^{(\eta_0)}_{\text{\rm Iw},h}(N,T)$ and define a meromorphic
$p$-adic function
$$
L_{\text{\rm Iw},h}(T,N,s)=f(\chi (\gamma_1)^s-1),
$$
where $\chi \,\,:\,\,\Gamma @>>>\Bbb Z_p^*$ is the cyclotomic character.
Let $\omega_N$ be a generator  of ${\det}_{\Bbb Z_p}(N).$ The isomorphism $D\simeq t_V(\Bbb Q_p)$ allows us to consider 
$\omega_N$ as a basis of ${\det}_{\Bbb Q_p} t_V(\Bbb Q_p).$ We also fix a generator
$\omega_T \in {\det}_{\Bbb Z_p}T^+$ and define the $p$-adic period $\Omega_p (\omega_T,\omega_{\text{\rm B}})\in \Bbb Q_p$ 
by 
$
\omega_{\text{\rm B}}=\Omega_p (\omega_T,\omega_{\text{\rm B}})\omega_T.
$
Our main result can be stated as follows.
\proclaim{Theorem 0.3} Let $V$ be a pseudo-geometric $p$-adic representation which is crystalline
at $p$. Assume that it  satisfies  conditions  1-4).  
Let $D$ be an admissible subspace of $\Dc (V).$  If  $\Cal L(D,V)\ne 0$ then

i)  $L_{\text{\rm Iw},h}(T,N,s)$ is a meromorphic $p$-adic function which has a zero 
at $s=0$ of order $e=\dim_{\Bbb Q_p}(D^{\Ph=p^{-1}}).$

ii) Let  $L_{\text{\rm Iw},h}^*(T,N,0)= \lim_{s\to 0} s^{-e}L_{\text{\rm Iw},h}(T,N,s)$
be the special value of $L_{\text{\rm Iw},h}(T,N,s)$ at $s=0.$ Then
$$
\multline
L_{\text{\rm Iw},h}^*(T,N,0)\overset{p}\to\sim \\
\Gamma (h)^{d_{+}(V)} \,\Cal L (D,V)\, 
\,E_p^*(V,1)\,{\det}_{\Bbb Q_p} \left (\frac{1-p^{-1}\Ph^{-1}}{1-\Ph} \,\vert D_{-1} \right ) 
\Omega_p(\omega_T,\omega_{\text{\rm B}})
i_{\omega_N,\omega_{\text{\rm B}}, p}\,(\Delta_{\text{\rm EP}} (T)),
\endmultline
$$
where $\Gamma (h)=(h-1)!$,  $d^+(V)=\dim_{\Bbb Q_p}(V^+)$,  $E_p(V,t)={\det} (1-\Ph t\mid \Dc (V))$ is the Euler factor at $p$ and   
$E_p^*(V,t)\,=\,E_p(V,t)\,\left (1-\dsize\frac{t}{p}\right )^{-e}.$
\endproclaim

\flushpar
{\bf Remarks 0.4.} 1) Assume that $V$ is an arbitrary pseudo-geometric representation
which is crystalline at $p$ and such that $\Dc (V)^{\Ph=1}=\Dc (V)^{\Ph=p^{-1}}=0.$  
In this case  the $p$-adic $L$-function has no trivial zeros (if exists) and  
a very general Iwasawa-theoretic descent result is proved in \cite{PR2}, Chapitre III.
If $V$ satisfies 1-4) and $\Dc (V)^{\Ph=p^{-1}}=0,$ it is easy to see that $\Cal L(D,V)=1$ and   Theorem 0.3 is a
particular case of this result, but our goal here is to study the case of trivial zeros.  
\flushpar
2) Let $E/\Bbb Q$ be an elliptic curve having  good reduction at $p.$ 
Consider the $p$-adic representation  $V=\text{\rm Sym}^2 (T_p(E))\otimes \Bbb Q_p$, where $T_p(E)$ is the $p$-adic Tate module 
of $E.$ It is easy to see that  $D=\Dc (V)^{\Ph=p^{-1}}$ is  one dimensional. In  this case some versions of Theorem 0.3 
were proved in \cite{PR3} and \cite{D} with  an ad hoc definition  of the $\Cal L$-invariant.
Remark that  $p$-adic $L$-functions attached to the symmetric square
of a newform were constructed by Dabrowski and Delbourgo \cite{DD}.
\flushpar
3) This theorem suggests that one should exist  an analytic $p$-adic $L$-function 
$L_{\text{\rm an}}(T,N,s)$ such that
$$
\align
&\bullet \,\, L_{\text{\rm an}}(T,N,s)\,\, \text{ has a zero of order $e-d^+(V)$  at $s=0$;}\\ 
&
\bullet \,\,L_{\text{\rm an}}^*(T,N,0)\overset{p}\to\sim  \,\Cal L (D,V)\, 
\,E_p^*(V,1)\,{\det}_{\Bbb Q_p} \left (\frac{1-p^{-1}\Ph^{-1}}{1-\Ph} \,\vert D_{-1} \right ) 
\frac{\Omega_p(\omega_T,\omega_{\text{\rm B}})}{ \Omega_\infty(\omega_N,\omega_B)}\,L(M,0).
\endalign
$$
\flushpar
{\bf 0.5.} The organization of the paper is as follows. In \S1 we review the theory of $(\Ph,\Gamma)$-modules,
in particular, the computation of cohomology of $(\Ph,\Gamma)$-modules of rank $1$ following \cite{C4}.
In \S2 we recall preliminaries on the Bloch-Kato exponential map and review the construction
of the large exponential map of Perrin-Riou given by Berger \cite{Ber3}. 
In \S3 we review the definition of the $\Cal L$-invariant given in \cite{Ben2} and interpret it
in terms of the Bockstein homomorphism associated  to the large exponential map.  In \S4 we prove 
Theorem 0.3 using the main result of \S3 and Nekov\'a\v r's Iwasawa-theoretic  descent techniques. 
In Appendix we prove derived versions of the well known computation of the local Galois cohomology
in terms of $(\Ph,\Gamma)$-modules  \cite{H1}, \cite{CC2}.
\newline
\newline  
{\bf Acknowledgements.} I am very grateful to  Jan Nekov\'a\v r and Daniel Delbourgo for several interesting discussions
and comments concerning this work. The main result of this paper was announced in a talk at the conference 
"Iwasawa 2008" (Ihrsee/Augsburg) organised by C. Greither and J. Ritter. I would like to thank them 
very much.

\head {\bf \S1. Preliminaries}
\endhead

{\bf 1.1. $(\Ph,\Gamma)$-modules.}

{\bf 1.1.1. The Robba ring} (see \cite{Ber1},\cite{C3}).
In this section  $K$ is a finite unramified extension of $\Bbb Q_p$ with residue field $k_K$, $O_K$ its ring of integers,
and $\sigma$ the absolute Frobenius of $K$. Let 
$\overline K$ an algebraic closure of $K$, $G_K=\text{Gal}(\bar K/K)$ and  $C$ the completion of $\overline K .$ 
Let  $v_p\,\,:\,\,C@>>>\Bbb R\cup\{\infty\}$ denote the $p$-adic valuation normalized so that $v_p(p)=1$ and 
set $\vert x\vert_p=\left (\frac{1}{p}\right )^{v_p(x)}.$ Write $B(r,1)$ for the $p$-adic annulus 
$B(r,1)=\{ x\in C \,\mid \, r\leqslant \vert x\vert <1 \}.$
As usually,  $\mu_{p^n}$  denotes the group of $p^n$-th roots of  unity.
Fix a system of primitive roots of unity   $\ep=(\zeta_{p^n})_{n\geqslant 0}$,
$\,\zeta_{p^n} \in \mu_{p^n} $ such that $\zeta_{p^n}^p=\zeta_{p^{n-1}}$ for all $n$.
Set $K_n=K(\zeta_{p^n})$, $K_{\infty}= \bigcup_{n=0}^{\infty}K_n$, $H_K=\Gal (\bar K/K_\infty)$, $\Gam =\G(K_{\infty}/K)$
and denote by $\chi \,:\,\Gam @>>>\Bbb Z_p^*$ the cyclotomic character.

Set 
$$
\widetilde {\bold E}^+=\varprojlim_{x\mapsto x^p} O_C/\,p\,O_C\,=\,\{x=(x_0,x_1,\ldots ,x_n,\ldots )\,\mid \,
x_i^p=x_i \,\,\forall i\in \Bbb N\}.
$$
Let $\hat x_n\in O_C$ be a lifting of $x_n$.
Then for all $m\geqslant 0$ the sequence $\hat x_{m+n}^{p^n}$ converges to 
$x^{(m)}=\lim_{n\to \infty} \hat x_{m+n}^{p^n}\in O_C$
which does not depend on the choice of  liftings.
 The ring $\widetilde {\bold E}^+$ equipped with the valuation $v_{\bold E}(x)=v_p(x^{(0)})$ is a
complete local ring of characteristic $p$ with  residue field
 $\bar k_K$. Moreover it is integrally closed in his field
of fractions $\widetilde {\bold E}=\text {\rm Fr}(\widetilde {\bold E}^+)$. 

Let $\widetilde \A=W(\widetilde \E)$ be the ring of Witt vectors with coefficients
in $\widetilde \E$. Denote by $[\,\,]\,:\,\widetilde \E@>>>W(\widetilde \E)$  the Teichmuller lift.
Any $u=(u_0,u_1,\ldots )\in \widetilde \A$ can be written in the form
$$
u=\underset{n=0}\to{\overset{\infty}\to \sum} [u^{p^{-n}}]p^n.
$$

Set $\pi=[\ep]-1$, $\A_{K_0}^+=O_{K_0} [[\pi]]$ and denote by $\A_{K}$  the $p$-adic completion
of $\A_{K}^+\left [1/{\pi}\right ]$. 
Let $\widetilde\boB= \widetilde \A\left [{1}/{p}\right ]$,  $\boB_{K}=\A_{K}\left [{1}/{p}\right ]$ and
let $\boB$ denote  the completion of the maximal unramified extension of $\boB_{K}$ in $\widetilde \boB$.
Set $\A=\boB\cap \widetilde \A$, $\widetilde \A^+=W(\E^+)$, $\A^+= \widetilde \A^+\cap \A$ and $\boB^+=\A^+\left [{1}/{p}\right ].$
All these rings are endowed with natural  actions of the Galois group $G_K$ and   Frobenius $\Ph$.

Set $\A_K=\A^{H_K}$ and  $\boB_K=\A_K\left [{1}/{p}\right ].$
Remark that $\Gamma$ and $\Ph$ act on $\boB_{K}$ by
$$
\aligned
& \tau (\pi)=(1+\pi)^{\chi (\tau)}-1,\qquad \tau \in \Gamma\\
&\Ph (\pi)=(1+\pi)^p-1.
\endaligned
$$
For any $r>0$ define
$$
\widetilde {\bold B}^{\dagger,r}\,=\,\left \{ x\in \widetilde
{\bold B}\,\,|\,\, \lim_{k\to +\infty} \left (
v_{\E}(x_k)\,+\,\dsize \frac{pr}{p-1}\,k\right )\,=\,+\infty
\right \}.
$$
Set ${\bold B}^{\dagger,r}=\boB \cap \widetilde\boB^{\dagger,r}$, 
$\boB_{K}^{\dagger,r}=\boB_{K} \cap \boB^{\dagger,r}$, 
${\bold B}^{\dagger}=\underset{r>0}\to \cup \boB^{\dagger,r}$
${\bold A}^{\dagger}=\bold A \cap {\bold B}^{\dagger}$
and $\bold B^\dag_K=\underset{r>0}\to \cup \boB_{K}^{\dagger,r}$.

It can be shown that for any $r\geqslant  p-1$

$$
\boB_{K}^{\dagger,r}=\left \{ f(\pi)=\sum_{k\in \Bbb Z}
a_k\pi^k\,\mid \, \text{\rm $a_k\in K$ and $f$ is holomorphic and bounded on $B(r,1)$} \right \}.
$$
Define
$$
\boB^{\dag,r}_{\text{rig},K}\,=\,\left \{ f(\pi)=\sum_{k\in \Bbb Z}
a_k\pi^k\,\mid \, \text{\rm $a_k\in K$ and $f$ is holomorphic  on $B(r,1)$} \right \}.
$$ 
Set    
$\CR (K) =\underset{r\geqslant  p-1}\to \cup \boB^{\dag,r}_{\text{rig},K}$ and $\CR^+(K)=\CR (K) \cap K[[\pi]].$ 
It is not difficult to check that these rings  are stable under $\Gamma$ and  $\Ph .$
To simplify notations we will write  $\Cal R=\CR (\Bbb Q_p)$ and $\Cal R^+= \CR^+(\Bbb Q_p).$
\newline
\,

{\bf 1.1.2. $(\Ph,\Gamma)$-modules} (see \cite{F2}, \cite{CC1}). 
Let $A$ be either $\boB_K^{\dag}$ or $\CR (K).$
A $(\Ph,\Gamma)$-module over A  is a finitely generated free $A$-module $D$
equipped with semilinear actions of $\Ph$ and  $\Gamma$ commuting to each other and such that
the induced linear map $\Ph \,:\,A\otimes_{\Ph} D @>>>D$ is an isomorphism.
Such a module    is said to be etale   if it admits a $\bold A_K^{\dag}$-lattice $N$ stable under $\Ph$ and $\Gamma$ and such that 
$\Ph \,:\,\bold A_K^{\dag}\otimes_{\Ph} N @>>>N$ is an isomorphism. The functor $D\mapsto \CR(K)\otimes_{\boB_K^{\dag}} D$
induces an equivalence between the category of etale $(\Ph,\Gamma)$-modules over $\boB_K^{\dag}$ and the category
of $(\Ph,\Gamma)$-modules over $\CR (K)$ which are of slope $0$ in the sense of Kedlaya's theory  (\cite{Ke} and
\cite{C5}, Corollary 1.5). Then Fontaine's classification of $p$-adic representations \cite{F2} together
with the main result of \cite{CC1} lead to the following statement.   

\proclaim{ Proposition 1.1.3} i) The functor 
$$
\bD^{\dagger}\,\,:\,\,V \mapsto \bD^{\dagger}(V)=(\boB^{\dagger}\otimes_{\Bbb Q_p}V)^{H_K}
$$
establishes an equivalence between the category of $p$-adic representations
of $G_K$ and the category of etale $(\Ph,\Gamma)$-modules over $\boB^{\dagger}_K .$

ii) The functor $\Ddagrig (V)=\Brigdag \otimes
_{\boB^{\dagger}_K} \bD^{\dagger}(V)$ 
gives  an equivalence between the
category of $p$-adic representations of $G_K$ and the category of
$(\Ph,\Gamma)$-modules over $\Brigdag$ of slope $0$.
\endproclaim
\demo{Proof} see \cite{C4}, Proposition 1.7.
\enddemo

{\bf 1.1.4.  Cohomology of $(\Ph,\Gamma)$-modules} (see \cite{H1}, \cite{H2}, \cite{Li}).
 Fix a generator $\gamma$ of $\Gamma$. If $D$ is a $(\Ph,\Gamma)$-module over $A$, we denote by
$ C_{\Ph,\gamma}(D)$ the complex
$$
C_{\Ph,\gamma} (D)\,\,:\,\,0@>f>>D @>>> D\oplus D@>g>> D@>>>0
$$
where $f(x)=((\Ph-1)\,x,(\gamma -1)\,x)$ and
$g(y,z)=(\gamma-1)\,y-(\Ph-1)\,z.$ Set $H^i(D)=H^i(C_{\Ph,\gamma}(D)).$ A short
exact sequence of $(\Ph,\Gamma)$-modules
$$
0@>>>D'@>>>D@>>>D''@>>>0
$$
gives rise to an exact cohomology sequence:
$$
0@>>>H^0(D')@>>>H^0(D)@>>>H^0(D'')@>>>H^1(D')@>>>\cdots @>>>
H^2(D'')@>>>0.
$$

\proclaim{Proposition 1.1.5} Let $V$ be a $p$-adic representation of $G_K.$
Then 

i) The complexes $\RG (K,V)$,  $C_{\Ph,\gamma}(\bD^{\dag}(V))$ and $C_{\Ph,\gamma}(\Ddagrig (V))$ are isomorphic in the derived
category of $\Bbb Q_p$-vector spaces $\Cal D(\Bbb Q_p).$

\endproclaim
\demo{Proof} This is a derived version of Herr's  computation of Galois cohomology \cite{H1}. 
The  proof is given in the Appendix, Propositions A.3 and  Corollary A.4.
\enddemo

{\bf 1.1.6.}  Recall that $\Lambda$ denotes the Iwasawa algebra of $\Gamma_1$, $\Delta=\Gal (K_1/K)$
and   $\Lambda (\Gamma)=\Bbb Z_p[\Delta]\otimes_{\Bbb Z_p}\Lambda$. Let $\iota \,\,:\,\,
\Lambda (\Gamma) @>>>\Lambda (\Gamma)$ denote the involution defined by $\iota (g)=g^{-1},$
$g\in \Gamma .$ 
If $T$ is a $\Bbb Z_p$-adic representation of $G_K$,  then the induced module
$\text{\rm Ind}_{K_\infty/K}(T)$ is isomorphic to $(\Lambda (\Gamma)\otimes_{\Bbb Z_p}T)^\iota$
and we set
$$
\RG_{\Iw}(K,T)\,=\,\RG (K, \text{\rm Ind}_{K_\infty/K}(T)).
$$
Write  $H^i_{\Iw}(K,T)$  for  the Iwasawa cohomology
$$
H^i_{\Iw}(K,T)=\varprojlim_{\text{cor}_{K_n/K_{n-1}}} H^i(K_n,T).
$$

Recall that there are canonical and functorial isomorphisms
 $$
 \aligned
 &\bold R^i \Gamma_{\Iw}(K,T)\,\simeq H^i_{\Iw}(K,T),\qquad i\geqslant 0,\\
 &\RG_{\Iw}(K,T)\otimes^{\bold L}_{\Lambda (\Gamma)} \Bbb Z_p[G_n] \simeq \RG (K_n,T)
 \endaligned
 $$
 (see  \cite{N2}, Proposition 8.4.22). The interpretation of the Iwasawa cohomology 
in terms of $(\Ph,\Gamma)$-modules was found by Fontaine (unpublished but see \cite{CC2}). 
We give here the derived version of this result. Let 
$\psi \,\,:\,\,\boB@>>>\boB$ be the operator
defined by the formula
$
\psi (x)\,=\,\frac{1}{p} \Ph^{-1} \left (\text{\rm Tr}_{\boB/\Ph (\boB)} (x)\right ).
$
We see immediately that  $\psi \circ  \Ph={\text{id}}.$ Moreover $\psi$ commutes with the action of $G_K$ 
and $\psi ({\bold A}^{\dag})={\bold A}^{\dag}$. Consider the complexes
$$
\aligned
&C_{\Iw,\psi}(T)\,\,:\,\, \bD (T)@>\psi-1>>\bD (T),\\
&C^{\dag}_{\Iw,\psi}(T)\,\,:\,\, \bD^{\dag} (T)@>\psi-1>>\bD^{\dag} (T).
\endaligned
$$

\proclaim{Proposition 1.1.7} i) The complexes $\RG_{\Iw}(K,T)$, 
$C_{\Iw,\psi}(T)$ and $C_{\Iw,\psi}^{\dag}(T)$ are naturally isomorphic in the derived category  $\Cal D(\Lambda (\Gamma))$
of $\Lambda (\Gamma)$-modules.
\endproclaim
\demo{Proof} See Proposition A.7 and Corollary A.8.
\enddemo

{\bf 1.1.8.} Finally, recall the computation of the cohomology of
$(\Ph,\Gamma)$-modules of rank $1$ following Colmez \cite{C4}. As in \cite{C4},  we
consider the case  $K=\Bbb Q_p$ and 
put  $\Cal R=\boB^{\dag}_{\text{\rm rig},\Bbb Q_p}$ and $\Cal R^+=\boB^+_{\text{\rm rig},\Bbb Q_p}.$
The differential operator
$\partial =(1+\pi)\dsize\frac{d}{d\pi}$ acts on $\CR$ and $\CR^+.$
If $\delta \,:\,\Bbb Q_p^*@>>> \Bbb Q_p^*$ is a continuous character, we write  $\Cal R (\delta)$
for  the $(\Ph,\Gamma)$-module $\Cal R e_{\delta}$ defined by
$\Ph (e_{\delta})=\delta (p) e_{\delta}$ and 
$\gamma (e_{\delta})=\delta (\chi(\tau))\,e_{\delta}.$ Let $ x$ denote
the character induced by the natural inclusion  of $\Bbb Q_p$ in $L$
and $\vert x\vert$ the character defined by $\vert
x\vert=p^{-v_p(x)}.$ 

\proclaim{Proposition 1.1.9} Let $\delta \,:\,\Bbb Q_p^*@>>>\Bbb Q_p^*$
be a continuous character. Then:

i) 

 $$H^0(\Cal R(\delta))=\cases \Bbb Q_p t^m &\text{\rm if
$\delta=x^{-m}$, $k\in \Bbb N$}\\
0 & \text{\rm otherwise.}
\endcases
$$

ii) 

$$\text{\rm dim}_{\Bbb Q_p} (H^1(\Cal R(\delta)))\,=\,\cases
2 &\text{\rm if either $\delta (x)= x^{-m}$, $m \geqslant 0$ or
$\delta (x)=\vert x\vert x^{m},$ $k\geqslant 1,$}\\
 1 &\text{\rm otherwise}.
\endcases
$$

iii) Assume that $\delta (x)= x^{-m},$ $m\geqslant 0.$ The classes 
$\cl (t^m,0)\,e_{\delta}$ and $\cl (0,t^m)\,e_{\delta}$ form a basis of  $H^1(\Cal R(x^{-m}))$.

iv) Assume that $\delta (x)= \vert x\vert x^{m},$ $m\geqslant 1.$ Then 
$H^1(\CR (\vert x\vert x^m)),$ $m \geqslant 1$  is generated by $\cl (\alpha_m)$ and $\cl (\beta_m)$
where
$$
\align
&\alpha_m=\frac{(-1)^{m-1}}{(m-1)!}\, \partial^{m-1} \left (\frac{1}{\pi}+\frac{1}{2},a \right )\,e_{\delta} ,\qquad
(1-\Ph)\,a=(1-\chi (\gamma) \gamma)\,\left (\frac{1}{\pi}+\frac{1}{2} \right ),\\
&\beta_m=\frac{(-1)^{m-1}}{(m-1)!}\,\partial^{m-1} \left (b,\frac{1}{\pi} \right )\,e_{\delta}, \qquad
(1-\Ph)\,\left (\frac{1}{\pi}\right )\,=\,(1-\chi (\gamma)\,\gamma)\,b
\endalign
$$
\endproclaim
\demo{Proof} See  \cite{C4}, sections 2.3-2.5.
\enddemo
\flushpar
\newline

{\bf 1.2. Crystalline representations.}

{\bf 1.2.1. The rings $\Bc$ and $\Bd$} (see \cite{F1}, \cite{F4}).
Let $\theta_0\,:\,\A^+@>>>O_C$ be the map given by the formula
$$
\theta_0 \left (\underset{n=0}\to{\overset{\infty}\to \sum} [u_n]p^n \right )\,=\,
\underset{n=0}\to{\overset{\infty}\to \sum} u_n^{(0)} p^n.
$$
It can be shown that $\theta_0$ is a surjective ring homomorphism and that $\ker (\theta_0)$ 
is the principal ideal generated by $\omega \,=\,\underset{i=0}\to {\overset{p-1} \to \sum} [\epsilon]^{i/p}.$
By linearity, $\theta_0$ can be extended to a map $\theta\,:\, \tilde {\boB}^+@>>>C$. The ring $\Bd^+$
is defined to be the completion of $\tilde {\boB}^+$ for the $\ker (\theta)$-adic topology:
$$
\Bd^+\,=\,\varprojlim_n \tilde {\boB}^+/ker (\theta)^n.
$$
This is a complete discrete valuation ring with residue field $C$ equipped with a natural
action of $G_K.$ Moreover, there exists a canonical embedding $\bar K  \subset \Bd^+.$
The series $t=\underset{n=0}\to{\overset \infty \to \sum}(-1)^{n-1}{\pi}^n/n$ converges
in the topology of $\Bd^+$ and it is easy to see that $t$ generates the maximal ideal of $\Bd^+.$ 
The Galois group acts on $t$ by the formula $g (t)=\chi (g)t.$  Let $\Bd=\Bd^+[t^{-1}]$ be the
field of fractions of $\Bd^+.$ This is a complete discrete valuation field equipped with 
a $G_K$-action and an exhaustive separated decreasing filtration $\F^i \Bd=t^i\Bd^+.$ 
As $G_K$-module, $\F^i\Bd/\F^{i+1}\Bd \simeq C(i)$ and  $\Bd^{G_K}=K.$

Consider the $PD$-envelope of $\A^+$ with a respect to the map $\theta_0$
$$
\A^{\text{\rm PD}}=\A^+\left [\frac{\omega^2}{2!},\frac{\omega^3}{3!},\ldots ,\frac{\omega^n}{n!},\ldots \right ]
$$
and denote by ${\A}^+_{\text{\rm cris}}$ its $p$-adic completion. Let $\Bc^+={\A}^+_{\text{\rm cris}}\otimes_{\Bbb Z_p}\Bbb Q_p$
and $\Bc=\Bc^+[t^{-1}].$ Then $\Bc$ is a subring of $\Bd$ endowed with  the induced filtration and  Galois action.
Moreover, it is equipped with a continuous Frobenius $\Ph$, extending the map $\Ph \,:\,\A^+@>>>\A^+.$
One has $\Ph (t)=p\,t.$ 
\newline
\, 

{\bf 1.2.2. Crystalline representations} (see \cite{F5}, \cite{Ber1}, \cite{Ber2}).
\newline
Let $L$ be a finite extension of $\Bbb Q_p$. Denote by $K$ its maximal
unramified subextension.  A filtered Dieudonn\'e module over $L$
is a finite dimensional $K$- vector space $M$ equipped with the
following structures:

$\bullet$ a $\sigma$-semilinear bijective map $\Ph \,:\,M@>>>M;$

$\bullet$ an exhaustive decreasing filtration $(\F^i M_L)$ on
the $L$-vector space $M_L=L\otimes_{K} M.$

A $K$-linear map $f\,:\,M@>>>M'$ is said to be a morphism of
filtered modules if

$\bullet$ $f(\Ph (d))=\Ph (f(d)),\qquad \text{\rm for all $d\in M;$}$

$\bullet$ $f(\F^i M_L)\subset \F^i M'_L, \qquad \text{\rm for
all $i\in \Bbb Z.$}$

The category  $\bold M\bold F^{\Ph}_L$ of filtered Dieudonn\'e modules is additive,
has kernels and cokernels but is not abelian. 
Denote by $\bold 1$  the vector space
$K_0$ with the natural action of $\sigma$ and  the
filtration given by
$$
\F^i\bold 1=\cases K,&\text{\rm if $i\leqslant 0,$}\\
0, &\text{\rm if $i>0.$}
\endcases
$$
Then $\bold 1$ is a unit object of $\bold M\bold F^{\Ph}_L$ i.e.
$
M\otimes \bold 1\simeq \bold 1\otimes M\simeq M
$
 for
any $M$.

If $M$ is a one dimensional Dieudonn\'e module and $d$ is a basis vector
of $M$, then $\Ph (d)=\alpha v$ for some $\alpha \in K$. Set
$t_N(M)= v_p(\alpha)$ and denote by $t_H(M)$ the unique filtration jump  
of $M.$ If $M$ is of an arbitrary finite dimension $d$, set 
$t_N(M)= t_N(\overset d\to \wedge M)$ and  $t_H(M)= t_H(\overset d\to \wedge M).$   
A Dieudonn\'e module $M$ is said to be weakly admissible if $t_H(M)=t_N(M)$
and if $t_H(M')\leqslant t_N(M')$  for any $\Ph$-submodule $M'\subset M$
equipped with the induced filtration. Weakly admissible modules 
form a subcategory of $\bold M\bold F_L$ which we denote by $\bold M\bold F_L^{\Ph,f}.$

If $V$ is a $p$-adic representation of $G_L$, define
$
\Dd (V)=(\Bd \otimes V)^{G_L}.
$
Then $\Dd (V)$ is a $L$-vector space  equipped with the  decreasing filtration
$\F^i\Dd (V)\,=\,(\F^i\Bd\otimes V)^{G_L}$.  One has 
$\dim_L \Dd (V)\leqslant \dim_{\Bbb Q_p}(V)$ and 
$V$ is said to be de
Rham if
$
\dim_L \Dd (V)= \dim_{\Bbb Q_p}(V).
$
Analogously one defines
$
\Dc (V)=(\Bc \otimes V)^{G_L}.
$
Then $\Dc (V)$ is a filtered Dieudonn\'e module over $L$ of
dimension $\dim_{K} \Dc(V)\leqslant \dim_{\Bbb Q_p}(V)$ and $V$ is
said to be crystalline if  the equality holds here.
In particular, for crystalline representations one has 
$
\Dd (V)= \Dc(V)\otimes_{K}L.
$
By the theorem of Colmez-Fontaine \cite{CF},  the functor $\Dc $
establishes an equivalence between the category of crystalline
representations of $G_L$ and $\bold M\bold F_L^{\Ph,f}$.
Its quasi-inverse $\bold V_{\text{\rm cris}}$ is given by
$
\bold V_{\text{\rm cris}}(D)=\F^0(D\otimes_{K_0} \Bc)^{\Ph=1}.
$

An important result of Berger (\cite{Ber 1}, Theorem 0.2)
says that
$\Dc (V)$ can be recovered from the $(\Ph,\Gamma)$-module $\Ddagrig (V).$
The situation is particularly simple  if 
If $L/\Bbb Q_p$ is unramified. In this case 
set 
${\bold D}^+(V)=(V\otimes_{\Bbb Q_p}\boB^+)^{H_K}$ and 
${\bold D}_{\text{\rm rig}}^+(V)= \CR^+(K)\otimes_{\boB_K^+}{\bold D}^+(V).$ Then
$$
\Dc (V)\,=\,\left ( {\bold D}_{\text{\rm rig}}^+(V) \left [\frac{1}{t}\right ]\right )^{\Gamma}
$$
(see \cite{Ber2}, Proposition 3.4).

\head {\bf \S 2. The  exponential map}
\endhead

{\bf 2.1. The Bloch-Kato exponential map} (\cite{BK}, \cite{N1}, \cite{FP}).

{\bf 2.1.1.}  Let $L$ be a finite extension of $\Bbb Q_p.$ Recall that we denote by $\bold M\bold F^{\Ph}_L$ 
the category of filtered Dieudonn\'e modules over $L.$ If $M$ is an object of $\bold M\bold F^{\Ph}_L$, define
$$
H^i (L,M)=\Ext^i_{\bold M\bold F^{\Ph}_L}(\bold 1, M), \qquad i=0,1.
$$
Remark that   $H^*(L,M)$ can be computed explicitly  as
the cohomology of the complex
$$
C^{\bullet}(M)\,\,:\,\,M@>f>>({M_L}/{\F^0 M_L}) \oplus M
$$
where the modules are placed in degrees $0$ and $1$ and $f(d)=
(d\pmod {\F^0 M_L},(1-\Ph)\,(d))$ (\cite{N1},\cite{FP}). Remark that if $M$ is weakly admissible
then each extension $0@>>>M@>>>M'@>>>\bold 1@>>>0$ is weakly admissible too and we can write
$ H^i (L,M)=\Ext^i_{\bold M\bold F^{\Ph,f}_L}(\bold 1, M).$ 
\newline
\,

{\bf 2.1.2.}  Let $\Rep_{\text{\rm cris}}(G_K)$ denote the category of crystalline representations
of $G_K.$ For any object $V$ of $\Rep_{\text{\rm cris}}(G_K)$ define
$$
H_f^i(K,V)=\Ext_{\Rep_{\text{\rm cris}}(G_K)}^i(\Bbb Q_p(0),V).
$$
An easy computation shows that
$$
H^i_f(K,V)=\cases H^0(K,V), &\text{if $i=0$,}\\
\ker\,(H^1(K,V)@>>>H^1(K,V\otimes \Bc)), &\text{if $i=1$,}\\
0, &\text{if $i\geqslant 2.$}
\endcases
$$
Let $t_V(K)=\Dd (V)/\F^0 \Dd (V)$ denote the tangent space of $V$. 
The rings $\Bd$ and $\Bc$ are related to each other via the fundamental
exact sequence
$$
0@>>>\Bbb Q_p@>>>\Bc @>f>>\Bd/\F^0\Bd \oplus \Bc @>>>0
$$
where $f(x)=( x\pmod{\F^0\Bd}, (1-\Ph)\,x)$ (see \cite{BK}, \S4).
Tensoring this sequence with $V$ and taking cohomology one obtains an exact sequence
$$
0@>>>H^0(K,V)@>>>\Dc (V)@>>>t_V(K)\oplus \Dc (V) @>>>H^1_f(K,V)@>>>0.
$$
The last map of this sequence gives rise to the  Bloch-Kato exponential map
$$
\exp_{V,K}\,\,:\,\, t_V(K)\oplus \Dc (V)@>>>H^1(K,V).
$$
Following \cite{F3}  set
$$
\RG_f(K,V)\,=\,C^\bullet (\Dc (V))\,=\,\left [ \Dc (V)@>f>> t_V(K)\oplus
\Dc (V) \right ].
$$

From the classification of crystalline representations in terms
of Dieudonn\'e modules  it follows that  the functor $\Vc$ induces natural isomorphisms
$$
r_{V,p}^i\,\,:\,\,\bold R^i\Gamma_f (K,V)@>>> H^i_f(K,V),\qquad i=0,1.
$$
The composite homomorphism
$$
t_K(V)\oplus \Dc (V)@>>> \bold R^1\Gamma_f (K,V)@>r_{V,p}^i>> H^i(K,V)
$$
 coincides with the Bloch-Kato exponential map
$\exp_{V,K}$ (\cite{N1}, Proposition 1.21).
\newline
\,

{\bf 2.1.3.} Let $g\,:\,B^{\bullet}@>>>C^{\bullet}$ be a morphism of complexes.
We denote by $\Tot^{\bullet} (g)$ the complex
$\Tot^n(g)=C^{n-1}\oplus B^{n}$ with differentials
$d^n\,:\,\Tot^n(g)@>>>\Tot^{n+1}(g)$ defined by the formula
\linebreak
$
d^n(c,b)=((-1)^n g^n(b)+d^{n-1}(c),d^n(b)).
$
It is well known that if
$
0@>>>A^{\bullet}@>f>>B^{\bullet}@>g>>C^{\bullet}@>>>0
$
is an exact sequence of complexes, then $f$
induces a quasi isomorphism
$
A^\bullet \overset{\sim}\to \rightarrow \Tot^{\bullet}(g).
$
In particular, tensoring the fundamental exact sequence with $V$, we obtain an
exact sequence of complexes
$$
0@>>>\RG (K,V)@>>>C_c^{\bullet} (G_K,V\otimes
\Bc)@>f>>C_c^{\bullet} (G_K,(V\otimes (\Bd /\F^0 \Bd)) \oplus
(V\otimes \Bc))@>>>0 
$$
which gives a quasi isomorphism  $ \RG (K,V) \overset{\sim}\to \rightarrow \Tot^{\bullet}
(f). $
Since $\RG_f(K,V)$ coincides tautologically with the complex
$$ C_c^0 (G_K,V\otimes \Bc)@>f>>C_c^0 (G_K,(V\otimes (\Bd /F^0 \Bd)) \oplus
(V\otimes \Bc))
$$
we obtain a diagram
$$
\xymatrix{ \RG (K,V) \ar[r]^{\sim} &\Tot^{\bullet} (f)   \\
&  \RG_f(K,V)
\ar[u]
\ar
@{.>}[ul]}
$$
 which defines a morphism 
$
\RG_f (K,V)@>>>\RG (K,V)
$
in $\Cal D(\Bbb Q_p)$ (see \cite{BF}, Proposition 1.17). Remark that the induced homomorphisms $\R^i\Gamma_f(K,V)@>>>
H^i(K,V)$ ($i=0,1$) coincide with the composition of $r_{V,p}^i$ 
with natural embeddings $H^i_f(K,V)@>>> H^i(K,V).$
\newline
\,

{\bf 2.1.4.} In this subsection we define an analogue of the exponential
map for crystalline $(\Ph,\Gamma)$-modules. Let  $K/\Bbb Q_p$ be an unramified extension. 
If $D$ is a $(\Ph,\Gamma)$-module over $\CR (K)$
define
$$
\Cal D_{\text{\rm cris}}(D)\,=\,\left (D\left [{1}/{t} \right ] \right )^{\Gamma}.
$$
It can be shown that   $\Cal D_{\text{\rm cris}}(D)$ is a finite dimensional $K$-vector space
equipped with a natural decreasing filtration $\F^i\Cal D_{\text{\rm cris}}(D)$ and
a semilinear action of $\Ph$. One says  that $D$ is crystalline if
$$
\dim_K(\Cal D_{\text{\rm cris}}(D))\,=\,\text{\rm rg}  (D).
$$
(see \cite{BC}).
From \cite{Ber4}, Th\'eor\`eme A it follows that the functor  $D\mapsto \Cal D_{\text{\rm cris}}(D)$ is an equivalence
between the category of crystalline $(\Ph,\Gamma)$-modules and $\bold M\bold F^{\Ph}_K.$
Remark that if $V$ is a $p$-adic representation of $G_K$ then 
$
\Dc (V)=\Cal D_{\text{\rm cris}}(\Ddagrig (V))
$ 
and $V$ is crystalline if and only if $\Ddagrig (V)$ is.

Let $D$ be a $(\Ph,\Gamma)$-module. To any cocycle $\alpha =(a,b)\in Z^1(C_{\Ph,\gamma}(D))$ 
one can associate the extension
$$
0@>>>D@>>>D_{\alpha}@>>> \CR(K)@>>>0
$$
defined by 
$$
D_{\alpha}=D\oplus \CR (K)\,e,\qquad (\Ph-1)\,e=a, \quad (\gamma-1)\,e=b.
$$
As usually, this gives rise to an isomorphism  $H^1(D)\simeq \text{\rm Ext}^1_{\CR} (\CR (K),D).$ 
We say that $\cl (\alpha)$ is crystalline if
$\dim_{K} \left (\Cal D_{\text{\rm cris}}(D_{\alpha})\right )=
\dim_{K} \left (\Cal D_{\text{\rm cris}}(D)\right )+1
$
and define 
$$
H^1_f(D)\,=\,\{ \cl (\alpha)\in H^1(D)\,\,\vert \,\, \text{\rm $\cl (\alpha)$ is crystalline}\,\}
$$
(see \cite{Ben2}, section 1.4.1). If $D$ is crystalline (or more generally potentially
semistable ) one has a natural isomorphism
$$
H^1(K, \Cal D_{\text{\rm cris}}(D))@>>>H^1_f(D).
$$
Set $t_D= \Cal D_{\text{\rm cris}}(D)/\F^0 \Cal D_{\text{\rm cris}}(D)$ and denote by
$
\exp_{D}\,\,:\,\,t_D\oplus \Cal D_{\text{\rm cris}}(D)@>>>H^1(D)
$
the composition of this isomorphism with the projection 
$t_D\oplus \Cal D_{\text{\rm cris}}(D) @>>>H^1(K, \Cal D_{\text{\rm cris}}(D))$ and the embedding
$H^1_f(D)\hookrightarrow H^1(D).$
\newline
\,

{\bf 2.1.5.} Assume  that $K=\Bbb Q_p.$ To simplify notation we will write $D_m$ for $\CR (\vert x\vert x^m)$ and $e_m$ for its canonical basis.
Then  $\Cal D_{\text{\rm cris}}(D_m)$ is the one dimensional $\Bbb Q_p$-vector space 
generated by $t^{-m} e_m$.   
As in \cite{Ben2}, we normalize the basis $(\cl (\alpha_m), \cl (\beta_m))$ of $H^1(D_m)$
putting $\alpha_m^*=\left (1-1/p\right )\,\cl (\alpha_m)$ and 
$\beta_m^*=\left (1-1/p \right ) \log (\chi (\gamma))\, \cl (\beta_m).$

\proclaim{Proposition 2.1.6} i) $H^1_f(D_m)$ is the one-dimensional $\Bbb Q_p$-vector space generated by 
$\alpha^*_m$. 
 
ii) The exponential map
$$
\exp_{D_m}\,\,:\,\,t_{D_m} @>>>H^1(D_m)
$$
sends $t^{-m}w_m$ to 
$- \,\alpha^*_m.
$
\endproclaim
\demo{Proof} This is a reformulation of \cite{Ben2}, Proposition 1.5.8 ii). 
\enddemo
$\,$

{\bf 2.2. The large exponential map.}

{\bf 2.2.1.} In this section $p$ is an odd prime number, $K$ is a finite unramified extension
of $\Bbb Q_p$ and $\sigma$ the absolute Frobenius acting on $K.$ 
Recall that $K_n=K(\zn)$ and $\dsize
K_\infty =\cup_{n=1}^{\infty} K_n.$ We set  $\Gamma =\Gal (K_\infty
/K),$  $\Gamma_n=\Gal (K_\infty/K_n)$ and $\Delta=\Gal (K_1/K)$ . 
Let  $\La=\Bbb Z_p[[\Gamma_1]]$ and $\Lambda (\Gamma)=\Bbb Z_p[\Delta]\otimes_{\Bbb Z_p}\Lambda.$
We will consider the following operators acting on the ring $K[[X]]$ of formal power series with coefficients in $K$:  

$\bullet$ The ring homomorphism $\sigma \,:\, K[[X]]@>>>K[[X]]$ defined
by
$
\sigma \left (\dsize \sum_{i=0}^\infty  a_iX^i\right )=\dsize\sum_{i=0}^\infty
\sigma (a_i)X^i;
$

$\bullet$  The ring homomorphism $\Ph \,:\,K[[X]]@>>>K[[X]]$ defined by
$$
\Ph \left (\sum_{i=0}^\infty  a_iX^i\right )=\sum_{i=0}^\infty
\sigma (a_i)\Ph(X)^i, \qquad \Ph (X)=(1+X)^p-1.
$$

$\bullet$ The differential operator $\dop=(1+X)\dsize\frac{d}{dX}.$ One
has
$
\dop\circ  \Ph=p\Ph \circ \dop.
$

$\bullet$ The operator $\psi \,: \,K[[X]]@>>> K[[X]]$ defined by
$
\psi (f(X))= \dsize\frac{1}{p} \Ph^{-1} \left
(\dsize\sum_{\z^p=1}f((1+X)\z-1)\right ).
$
It is easy to see that $\psi $ is a left inverse to $\Ph,$ i.e.
that
$
\psi \circ \Ph=\text{\rm id}.
$

$\bullet$  An action of $\Gamma$ given by
$
\g \left (\dsize\sum_{i=0}^\infty  a_iX^i\right )=\dsize \sum_{i=0}^\infty a_i
\g( X)^i, \qquad \g (X)=(1+X)^{\chi (\g)}-1.
$

Remark that these formulas are compatible with the definitions from sections 1.1.1 and 1.1.6.
Fix a generator $\g_1\in \Gamma_1$ and define 
$$
\Cal H\,=\,\{f(\g_1-1)\,|\,f \in \Bbb Q_p[[X]] \,\,
\text{{\rm is holomorphic on $B(0,1)$}}\},\qquad
\Cal H(\Gamma)\,=\,\Bbb Z_p[\Delta]\otimes_{\Bbb Z_p}\Cal H.
$$

{\bf 2.2.2.} It is well known that  $ \Bbb Z_p[[X]]^{\psi=0}$ is a free $\Lambda$-module generated by
$(1+X)$  and the operator $\dop$ is bijective on $\Bbb
Z_p[[X]]^{\psi=0}.$ If $V$ is  a crystalline representation of $G_K$ 
put
$
\Cal D(V)=\Dc (V)\otimes _{\Bbb Z_p} \Bbb Z_p[[X]]^{\psi=0} .
$
Let 
$\boldsymbol \Xi_{V,n}^{\ep}\,:\,\Cal D(V)_{\Gamma_n} [-1]@>>>\RG_f
(K_n,V)
$
be the map defined by
$$
 \boldsymbol{\Xi}_{V,n}^\ep (\alpha)\,=\,\cases p^{-n} (\sum_{k=1}^n
(\sigma \otimes \Ph)^{-k} \alpha (\zeta_{p^k}-1),\,-\alpha (0))
&\text{\rm if
$n\geqslant 1$,} \\
\text{\rm Tr}_{K_1/K}\,\left ( \boldsymbol{\Xi}_{V,1}^{\ep} (\alpha)\right )
&\text{\rm if $n=0.$}
\endcases
$$
An easy computation shows that $ \boldsymbol{\Xi}_{V,0}^\ep \,:\,
\Dc (V) [-1] @>>> \RG_f(K,V)$ is given by the formula
$$
\boldsymbol{\Xi}_{V,0}^\ep
(a)\,=\,\frac{1}{p}\,(-\Ph^{-1}(a),-(p-1)\,a).
$$
In particular, it is homotopic to the map
$
a\mapsto -(0,(1-p^{-1}\Ph^{-1} )\,a).
$
Write
$$
\Xi_{V,n}^\ep \,:\,\Cal D(V)@>>> \bold R^1 \Gamma
(K_n,V)=\frac{t_V(K_n)\oplus \Dc (V)}{\Dc (V)/V^{G_K}}
$$
denote the homomorphism  induced by $\boldsymbol \Xi_{V,n}^\ep.$ Then 
$$
\Xi_{V,0}^\ep (a)=-(0,(1-p^{-1}\Ph^{-1} )\,a) \pmod{\Dc (V)/V^{G_K}}.
$$
If $\Dc (V)^{\Ph=1}=0$ the operator $1-\Ph$ is invertible on $\Dc (V)$ and we can write  
$$
\Xi_{V,0}^\ep (a)= \left (\frac{1-p^{-1}\Ph^{-1}}{1-\Ph} \,a,0 \right ) \pmod{\Dc (V)/V^{G_K}}. \tag{2.1}
$$

\flushpar
For any $i\in \Bbb Z$ let
$\Delta_i\,:\, \Cal D(V)@>>>\dsize\frac{\Dc(V)}{(1-p^i\Ph)\Dc
(V)}\,\otimes \Bbb Q_p(i)
$
be the map given by 
$$\Delta_i(\alpha (X))= \dop^i\alpha(0)\otimes \ep^{\otimes i}\pmod
{(1-p^i\Ph)\Dc(V)}.$$
 Set $\Delta=\oplus_{i\in \Bbb Z}\Delta_i .$
  If $\alpha\in \Cal
D(V)^{\Delta =0},$ then by  \cite{PR1}, Proposition 2.2.1 there exists $F\in \Dc(V)\otimes_{\Bbb Q_p}\Bbb Q_p[[X]]$ 
which converges on the open unit disk 
and such
that $ (1-\Ph)F=\alpha.$ A short computation shows that
$$
\Xi_{V,n}^\ep (\alpha) = p^{-n} ( (\sigma\otimes \Ph)^{-n}
(F)(\zn-1),0)\pmod {\Dc (V)/V^{G_K}},\qquad \text{\rm
if $n\geqslant 1$}
$$
(see \cite{BB}, Lemme 4.9). 
\newline
\newline
{\bf 2.2.3.} As $\Bbb Z_p[[X]] \left [1/p\right
] $ is a principal ideal domain and $\Cal H$ is $\Bbb Z_p[[X]]
\left [1/p\right ]$-torsion free,  $\Cal
H$ is flat. Thus
$$
C_{\Iw, \psi}^\dag(V)\otimes^{\bold L}_{\Lambda_{\Bbb
Q_p}}\HG\,=\,C_{\Iw, \psi}^\dag(V)\otimes_{\Lambda_{\Bbb Q_p}
}\HG\,=\,\left [ \HG\otimes_{\Lambda_{\Bbb Q_p}} \bD^\dag
(V)@>\psi-1>>\HG\otimes_{\Lambda_{\Bbb Q_p}}\bD^\dag (V) \right ].
$$
By proposition 1.1.7 on has an isomorphism in $\Cal D (\HG)$
$$
\RG_{\Iw}(K,V)\otimes^{\bold L}_{\Lambda_{\Bbb Q_p}}\HG \simeq
C_{\Iw, \psi}^\dag (V)\otimes_{\Lambda_{\Bbb Q_p}}\HG .
$$
The action of $\HG$ on
$\bD^{\dag}(V)^{\psi=1}$ induces an injection
$
\HG \otimes_{\La_{\Bbb Q_p}} \bD^{\dag}(V)^{\psi=1}
\hookrightarrow \Ddagrig (V)^{\psi=1}.
$
Composing this map with the canonical isomorphism
$\Hi^1(K,V)\simeq \bD^{\dag}(V)^{\psi=1}$ we obtain a map
$
\HG \otimes_{\La_{\Bbb Q_p}} \Hi^1(K,V) \hookrightarrow \Ddagrig
(V)^{\psi=1}.
$
For any $k\in \Bbb Z$ set $\nabla_k=t\dop-k=
t\dsize\frac{d}{dt}-k.$ An easy induction shows that
$
\nabla_{k-1}\circ \nabla_{k-2}\circ \cdots \circ \nabla_0\,=\,t^k\dop^k.
$

Fix $h\geqslant 1$  such that $\F^{-h}\Dc (V)=\Dc (V)$ and $V(-h)^{G_K}=0.$
For any  $\alpha \in \Cal D(V)^{\Delta=0}$ define 
$$
\Omega_{V,h}^\ep(\alpha) \,=\,(-1)^{h-1}\frac{\log \chi (\gamma_1)}{p}\,\nabla_{h-1}\circ
\nabla_{h-2}\circ \cdots \nabla_0 (F (\pi)),
$$
where $F\in \Cal H (V)$ is such that $(1-\Ph)\,F= \alpha.$ 
It is easy to see that $\Omega_{V,h}^\ep(\alpha) \in
\Drig^+(V)^{\psi=1}.$ In \cite{Ber3} Berger shows that
$\Omega_{V,h}^\ep(\alpha)\in \HG \otimes_{\La_{\Bbb
Q_p}}\bD^{\dag} (V)^{\psi=1}$ and therefore  gives rise to a map
$$
\bExp_{V,h}^\ep\,:\,\Cal D(V)^{\Delta=0}[-1] @>>>
\RG_{\Iw}(K,V)\otimes^{\bold L}_{\La_{\Bbb Q_p}}\HG
$$
Let
$$
\Exp_{V,h}^\ep\,:\,\Cal D(V)^{\Delta=0}@>>>\HG\otimes_{\La_{\Bbb
Q_p}}H_{\Iw}^1(K,V)
$$
denote the map induced by $\bExp_{V,h}^\ep$ in degree $1$. The
following theorem is a reformulation of the construction of the large exponential
map given by Berger in \cite{Ber3}.

\proclaim{Theorem 2.2.4} Let
$$
\bExp_{V,h,n}^\ep \,\,:\,\,\Cal D(V)^{\Delta=0}_{\Gam_n}[-1] @>>>
\RG_{\Iw} (K,V)\otimes^{\bold L}_{\La_{\Bbb Q_p}}\Bbb Q_p[G_n].
$$
denote the map induced by $\bExp_{V,h}^\ep.$ Then for any
$n\geqslant 0$ the following diagram in $\Cal D (\Bbb Q_p[G_n])$
is commutative:

$$
\xymatrix{  \Cal D(V)^{\Delta=0}_{\Gamma_n}[-1] \ar[rr]^
{\bExp^{\ep}_{V,h,n}} \ar[d]^{\boldsymbol \Xi^{\ep}_{V,n}} & &
\RG_{\Iw} (K,V)\otimes^{\bold L}_{\La_{\Bbb Q_p}}\Bbb Q_p[G_n] \ar
[d]^{\simeq}
\\
\bold R\Gamma_f (K_n,V) \ar[rr]^{(h-1)!} & &\RG (K_n,V)\,. }
$$

In particular,  $\Exp_{V,h}^\ep$ coincides with the large
exponential map of Perrin-Riou.
\endproclaim
\demo{Proof} Passing to cohomology in the previous diagram one
obtains  the  diagram
$$
\CD \Cal D(V)^{\Delta=0} @> \Exp^{\ep}_{V,h}
 >> \Cal H (\Gam)\otimes_{\La_{\Bbb Q_p}}H^1_{\text{\rm Iw}}(K,V)
\\
@V\varXi^{\ep}_{V,n}VV         @VV\text {\rm pr}_{V,n}V
\\
\bD_{\text{dR}/K_n}(V)\oplus \Dc (V) @>(h-1)!\,\exp_{V,K_n}>>
H^1(K_n, V)
\endCD
$$
which is exactly the definition of the large exponential map. Its
commutativity is proved in \cite{Ber3}, Theorem II.13. Now, the
theorem is an immediate consequence of the following remark. Let
$D$ be a free $A$-module and let $f_1,f_2\,:\,D[-1] @>>>
K^\bullet$ be two maps from $D[-1]$ to a complex of $A$-modules
such that the induced maps $h_1(f_1)$ and $h(f_2)\,:\,
D@>>>H^1(K^\bullet)$ coincide. Then $f_1$ and $f_2$ are homotopic.
\enddemo

\head {\bf \S3. The $\Cal L$-invariant}
\endhead

{\bf 3.1. Definition of the $\Cal L$-invariant} (\cite{Ben2}). 

{\bf 3.1.1.} In this section we recall the definition of the $\Cal
L$-invariant for the case of crystalline
representations. For further details and proofs see \cite{Ben2}, \S2.
Let $S$ be a finite set of primes  of $\Bbb Q$
containing $p$  and  $G_S$  the Galois group of the maximal
algebraic extension of $\Bbb Q$ unramified outside $S\cup\{\infty\}$. For each place
$v$ we denote by $G_v$  the decomposition at $v$ group and
by $I_v$ and $f_v$ the inertia subgroup and Frobenius automorphism
respectively.
Let $V$ be a $p$-adic pseudo-geometric representation of $G_S$.  
Thus $V$ is a de Rham at $p$. For any $v\notin
\{p,\infty\}$  set
 $$
 \RG_f (\Bbb Q_v,V)\,=\,\left [ V^{I_v} @>1-f_v>> V^{I_v} \right
 ],
 $$
where the terms are placed in degrees $0$ and $1$ (see \cite{F3}, \cite{BF}).  Observe that
there is a natural quasi-isomorphism $\RG_f(\Bbb Q_v,V)\simeq
C_c^\bullet (G_v/I_v, V^{I_v}).$ In particular, $\bold R^0\Gamma
(\Bbb Q_v,V)=H^0(\Bbb Q_v,V)$ and $\bold R^1\Gamma_f (\Bbb
Q_v,V)=H^1_f(\Bbb Q_v,V)$ where
$$
H^1_f(\Bbb Q_v,V)=\ker (H^1(\Bbb Q_v,V)@>>>H^1(\Bbb Q_v^{\text{\rm
ur}},V)).
$$
For $v= p$ the complex $\RG_f (\Bbb Q_v,V)$ was defined in section 2.1.2.
To simplify notation write $H^i_S(V)=H^i(G_S,V).$ The Selmer
group of $V$ is defined by
$$
H^1_f(V)\,=\,\ker \left ( H^1_S(V)@>>>\bigoplus_{v\in S}
\frac{H^1(\Bbb Q_v,V)}{H^1_f(\Bbb Q_v,V)} \right ).
$$

{\bf 3.1.2.} Assume that $V$ satisfies the following conditions:
\newline
\,

{\bf C1)} $H^1_f(V)=H^1_f(V^*(1))=0$.

{\bf C2)}  $H^0_S(V)=H^0_S(V^*(1))=0.$

{\bf C3)} $V$ is crystalline at $p$, $\Dc (V)^{\Ph=1}=0$  and the linear
map $1-p^{-1} \Ph^{-1} \,:\,\Dc (V)@>>>\Dc
(V)$ is semisimple.

{\bf C4)} The $(\Ph,\Gamma)$-module $\Ddagrig (V)$ has no
crystalline subquotient of the form
$$
0@>>>\Cal R (\mid x\mid x^k)@>>>U@>>>\Cal R@>>>0, \qquad
k\geqslant 1.
$$

Write $c$ for the complex conjugation and set  $d_{\pm} (V)=\dim (V^{c=\pm 1}).$ 
From the Poitou-Tate exact sequence it follows that $\dim_{\Bbb Q_p}t_V(\Bbb Q_p)\,=\,d_+(V).$ 
We say that a $\Bbb Q_p$-subspace $D\subset
\Dc (V)$ is admissible if it is stable under  $\Ph$
and the natural projection $D@>>>t_V(\Bbb Q_p)$ is an isomorphism.
\newline
\,

{\bf 3.1.3.} Let $D$ be an admissible subspace of $\Dc (V).$ As $1-p^{-1}\Ph^{-1}$ acts semisimply,
one has a decomposition
$
D\simeq D_{-1} \oplus D^{\Ph=p^{-1}}
$
where $D_{-1}=(\Ph-p^{-1})\,D $ is stable under $\Ph $ and
$(D_{-1})^{\Ph=p^{-1}}=0.$ Consider the filtration $(D_i)$ on $\Dc (V)$ defined by
$$
D_i=\cases 0&\text{if $i=-2$},\\
D_{-1} &\text{if $i=-1$},\\
D &\text{if $i=0$},\\
\Dc (V) &\text{if $i=1$}.
\endcases
$$
By Berger's theory \cite{Ber4} $(D_i)$ induces a filtration on $\Ddagrig (V)$:
$$
0\subset F_{-1}\Ddagrig (V)\subset F_0\Ddagrig (V)\subset F_1\Ddagrig
(V)=\Ddagrig (V).
$$
Explicitely
$ F_{i}\Ddagrig (V)\,=\,\Ddagrig (V)\cap \left (
D_{i}\otimes_{\Bbb Q_p} \Cal R\left [{1}/t \right ] \right )
$
(\cite{BC}, section 2.4.2).
Set  ${\text{\rm gr}}_i \Ddagrig (V)\,=\, F_i\Ddagrig (V)/F_{i-1} \Ddagrig
(V)$. By \cite{Ben2}, Corollary 1.4.6 the exact sequence 
$$
0@>>> F_0 \Ddagrig (V) @>>>  \Ddagrig (V) @>>> {\text{\rm gr}}_1
\Ddagrig (V) @>>>0
$$
gives rise to  exact sequences
$$
\cdots @>>> H^0({\text{\rm gr}}_1 \Ddagrig (V))@>>> H^1(F_0
\Ddagrig (V)) @>>> H^1 (\Ddagrig (V)) @>>>H^1({\text{\rm gr}}_1 \Ddagrig (V))@>>> \cdots 
$$
and
$$
\cdots @>>> H^0({\text{\rm gr}}_1 \Ddagrig (V))@>>> H^1_f(F_0
\Ddagrig (V)) @>>> H^1_f (\Ddagrig (V)) @>>> H^1_f({\text{\rm
gr}}_1 \Ddagrig (V))@>>>0
$$
The condition {\bf C3)} implies that $\Cal D_{\text{\rm cris}}
({\text{\rm gr}}_1 \Ddagrig (V))^{\Ph=1}=0.$ Since $D$ is
admissible, the Hodge-Tate weights of ${\text{\rm gr}}_1 \Ddagrig
(V)$ are $\leqslant 0$ and  by Proposition 1.4.4 of \cite{Ben2} 
 $H^0({\text{\rm gr}}_1 \Ddagrig (V))=0$ and
$H^1_f({\text{\rm gr}}_1 \Ddagrig (V))=0.$ This shows that
$H^1(F_0 \Ddagrig (V))$ injects into $H^1 (\Ddagrig
(V))\simeq H^1(\Bbb Q_p,V)$ and that 
$$
H^1_f(F_0 \Ddagrig (V)) \simeq H^1_f (\Ddagrig (V)) \simeq
H^1_f(\Bbb Q_p,V).
$$

Now, consider the short exact sequence
$$
0@>>> F_{-1} \Ddagrig (V) @>>>  F_0\Ddagrig (V) @>>> {\text{\rm
gr}}_0 \Ddagrig (V) @>>>0.
$$
Since $\Cal D_{\text{\rm cris}} (F_{-1} \Ddagrig (V))^
{\Ph=p^{-1}}=0$ and  Hodge-Tate weights of $F_{-1} \Ddagrig (V)$
are $>0$ we have 
$H^1_f(F_{-1}\Ddagrig (V))=H^1(F_{-1}\Ddagrig (V))$ by \cite{Ben2}, Proposition 1.4.4.
As $ \Cal D_{\text{\rm cris}}((F_{-1}\Ddagrig (V))^*(\chi))$ is dual to 
$\Cal D_{\text{\rm cris}} (F_{-1} \Ddagrig (V))$, the map $1-\Ph$ is bijective on  
 $ \Cal D_{\text{\rm cris}}((F_{-1}\Ddagrig (V))^*(\chi))$ and 
\linebreak
 $H^0((F_{-1}\Ddagrig (V))^*(\chi))=0.$ Using the local duality \cite{Li}
we obtain that $H^2(F_{-1}\Ddagrig (V))=0$. Finally 
$\Cal D_{\text{\rm cris}}{(\text{\rm gr}}_0 \Ddagrig (V))^{\Ph=1}=0$ implies that
$H^0(\text{\rm gr}_0\Ddagrig (V))=0$. Thus we have
exact sequences
$$
\align & 0@>>> H^1( F_{-1} \Ddagrig (V)) @>>> H^1( F_0\Ddagrig
(V)) @>>> H^1({\text{\rm gr}}_0 \Ddagrig (V)) @>>>0,\\
& 0@>>> H^1( F_{-1} \Ddagrig (V)) @>>> H^1_f( F_0\Ddagrig (V))
@>>> H^1_f({\text{\rm gr}}_0 \Ddagrig (V)) @>>>0.
\endalign
$$
Therefore
$$
H^1({\text{\rm gr}}_0 \Ddagrig (V))\simeq \frac{H^1( F_0\Ddagrig
(V))}{H^1( F_{-1} \Ddagrig (V))} \hookrightarrow \frac{H^1(\Bbb
Q_p,V)}{H^1( F_{-1} \Ddagrig (V))}
$$
and
$$
H^1_f({\text{\rm gr}}_0 \Ddagrig (V))\simeq \frac{H^1_f(\Bbb
Q_p,V)}{H^1( F_{-1} \Ddagrig (V))}.
$$

As $\Cal D_{\text{\rm cris}}{(\text{\rm gr}}_0 \Ddagrig (V))^{\Ph=p^{-1}}=\Cal D_{\text{\rm cris}}{(\text{\rm gr}}_0 \Ddagrig (V)),$
Proposition 1.5.9 of \cite{Ben2} implies that 
$$
\text{\rm gr}_0 \Ddagrig (V)\simeq \underset{i=1}\to{\overset
e\to\oplus} D_{m_i},  \qquad e=\dim_{\Bbb Q_p} (D^{\Ph=p^{-1}})
$$
where $D_{m_i} = \Cal R(\vert x\vert x^{m_i}),$ $m_i\geqslant 1.$ 
By Proposition 2.1.6  $H^1_f(D_m)$ is    generated
by $\alpha^*_m $ and we denote by $H^1_c(D_m)$ the subspace generated by $\beta^*_m.$
This gives  a  decomposition
$$
H^1(\text{\rm gr}_0 \Ddagrig (V))\simeq H^1_f(\text{\rm gr}_0 \Ddagrig (V))\oplus
H^1_c(\text{\rm gr}_0 \Ddagrig (V)).
$$
In particular,  $H^1_f(\text{\rm gr}_0 \Ddagrig (V))$ and   $H^1_c(\text{\rm gr}_0 \Ddagrig (V))$
are $\Bbb Q_p$-vector spaces of dimension $e.$
Further, fixing the basis $\alpha_m^*, \beta_m^*$ of $H^1(D_m)$  we fixe isomorphisms 
$$i_{D,f}\,:\,\Cal D_{\text{cris}}(\text{\rm gr}_0 \Ddagrig (V))   \iso H^1_f(\text{\rm gr}_0 \Ddagrig (V)),\qquad
i_{D,c}\,:\,\Cal D_{\text{cris}}(\text{\rm gr}_0 \Ddagrig (V)) \iso H^1_c(\text{\rm gr}_0 \Ddagrig (V)). 
$$
The condition {\bf C1)} together with the Poitou-Tate exact sequence implies that
$$
H^1_S(V) \simeq \bigoplus_{v\in S} \frac{H^1(\Bbb
Q_v,V)}{H^1_f(\Bbb Q_v,V)}.
$$
Let $H^1_S(D,V)$ be the subspace of $H^1_S(V)$
whose image under this isomorphism is 
\linebreak
$H^1( F_0\Ddagrig (V))/
 H^1_{f}(\Bbb Q_p,V).$ The localization
map $H^1_S(D,V)@>>>\dsize \frac{H^1(\Bbb Q_p,V)}{H^1(F_{-1}\Ddagrig
(V))}$ is injective and its image is contained in $H^1(\text{\rm
gr}_0 \Ddagrig (V)).$
Hence, we have a diagram
$$
\xymatrix{
\Cal D_{\text{\rm cris}}({\text{\rm gr}}_0 \Ddagrig (V)) \ar[r]^{\overset{i_{D,f}}\to \sim} & H^1_f({\text{\rm gr}}_0 \Ddagrig (V))\\
H^1_S(D,V) \ar[u]^{\rho_{D,f}} \ar[r] \ar[d]_{\rho_{D,c}} &
H^1({\text{\rm gr}}_0 \Ddagrig (V)) \ar[u]_{p_{D,f}}
\ar[d]^{p_{D,c}}
\\
\Cal D_{\text{\rm cris}}({\text{\rm gr}}_0 \Ddagrig (V)) \ar[r]^{\overset{i_{D,c}}\to\sim} &H^1_c({\text{\rm gr}}_0 \Ddagrig (V)),}
$$
where $\rho_{D,f}$ and $\rho_{D,c}$ are defined as the unique maps making
this diagram commute.
From the definition of $H^1_S(D,V)$ it follows that  $\rho_{D,c}$ is an isomorphism.

\proclaim{Definition 3.1.4} The determinant
$$
\Cal L (V,D)= \det \left ( \rho_{D,f} \circ \rho^{-1}_{D,c}\,\mid \,\Cal
D_{\text{\rm cris}}(\text{\rm gr}_0 \Ddagrig (V)) \right )
$$
will be called  the $\Cal L$-invariant associated to $V$ and $D$.
\endproclaim
\flushpar
{\bf 3.2. The Bockstein homomorphism.}
\newline
{\bf 3.2.1.} In this section we interpret  $\Cal L(D,V)$ in terms
of the Bockstein homomorphism associated to the large exponential map.
This interpretation is crucial for the  proof of the main theorem of this paper.
Recall that  $H^1(\Bbb Q_p,\Cal H(\Gamma)\otimes_{\Bbb Q_p}V)= \Cal H(\Gamma)\otimes_{\Lambda (\Gamma)} H_{\Iw}^1(\Bbb Q_p,V)$
injects into $\Ddagrig (V).$  Set 
$$F_iH^1(\Bbb Q_p,\Cal H(\Gamma)\otimes_{\Bbb Q_p}V)=  F_i\Ddagrig (V)
\cap H^1(\Bbb Q_p,\Cal H(\Gamma)\otimes_{\Bbb Q_p}V). $$ 
As in section 2.2 we fix a generator $\gamma \in \Gamma$. 

\proclaim{Proposition 3.2.2} Let $D$ be an admissible subspace of
$\Dc (V).$ For any  $a \in D^{\Ph=p^{-1}}$ let $\alpha  \in \Cal D(V)$
be such that $\alpha (0)=a .$ Then

i) There exists a unique $\beta \in F_0 H^1(\Bbb Q_p,\Cal
H(\Gamma)\otimes V)$ such that
$$
(\gamma-1)\,\beta =\Exp_{V,h}^\ep (\alpha ).
$$

ii) The composition map
$$
\align &\delta_{D,h}\,:\,D^{\Ph=p^{-1}}@>>>F_0 H^1(\Bbb Q_p,\Cal
H(\Gamma)\otimes V)@>>>H^1(\text{\rm gr}_0 \Ddagrig (V))\\
&\delta_{D,h} (a)\,=\,\beta  \pmod{H^1(F_{-1} \Ddagrig (V))}
\endalign
$$
is  given
explicitly by  the following formula: 
$$
\delta_{D,h}(\alpha)\,=\,-(h-1)!\,\left (1-\frac{1}{p}\right )^{
 -1}(\log \chi
(\gamma))^{-1}\, i_{D,c} (\alpha).
$$
\endproclaim
\demo{Proof}  Since $\Dc (V)^{\Ph=1}=0,$ the operator $1-\Ph$ is invertible on
$\Dc (V)$ and we have a diagram
$$
\xymatrix{ \Cal D(V)^{\Delta=0} \ar[rr]^{\Exp_{V,h}^\ep}
\ar[d]^{\Xi_{V,0}^\ep} & & H^1(\Bbb Q_p, \Cal H(\Gamma)\otimes V)
\ar[d]^{\text{\rm pr}_V}\\
\Dc (V) \ar[rr]^{(h-1)!\exp_V} & & H^1(\Bbb Q_p,V). }
$$
where $\Xi_{V,0}^\ep
(\alpha)=\dsize\frac{1-p^{-1}\Ph^{-1}}{1-\Ph}\,\alpha (0)$
(see (2.1)). If
$\alpha \in D^{\Ph=p^{-1}}\otimes \Bbb Z_p[[X]]^{\psi=0},$ then
$\Xi_{V,0}^\ep(\alpha)=0$ and 
$
\text{\rm pr}_V \left (\Exp_{V,h}^\ep (\alpha)\right )\,=\,0.
$
On the other hand, as $\Dc (V)^{\Ph=1}=0$, we have $V^{G_K}=0$ and the
 map 
$
\left (\Cal H(\Gamma)\otimes_{\Lambda_{\Bbb Q_p}} \Hi^1(\Bbb Q_p,V)\right )_{\Gamma} @>>> H^1(\Bbb Q_p,V)
$
is injective.  Thus there exists a unique 
\linebreak
$\beta \in \Cal H(\Gamma)\otimes_{\Lambda} \Hi^1 (\Bbb Q_p,T)$
such that 
$
\Exp_{V,h}^\ep (\alpha) \,=\,(\gamma-1)\,\beta.
$
Now take $a\in D^{\Ph=p^{-1}}$ and set
$$
f= a \otimes \ell \left (\frac{(1+X)^{\chi (\gamma)}-1}{X}\right ),
$$
where  
$
\ell (g)\,=\,\dsize\frac{1}{p}\log \left(\frac{g^p}{\Ph (g)}\right ). 
$
 An easy computation shows that
$$
\sum_{\zeta^p=1} \ell \left ( \frac{\zeta^{\chi (\gamma)}(1+X)^{\chi (\gamma)}-1}{\zeta (1+X)-1} \right )\,=\,0.
$$
Thus $f \in D^{\Ph=p^{-1}}\otimes \Bbb Z_p[[X]]^{\psi=0}.$
Write $\alpha$ in the form
$\alpha\,=\,(1-\Ph)\,(1-\gamma)\,(a\otimes \log (X)).$ Then 
$$
\Omega_{V,h}(\alpha)\,=\,(-1)^{h-1}\frac{\log \chi (\gamma_1)}{p}\,t^h\partial^h((\gamma-1)\,(a\log (\pi))\,=\,
\frac{\log \chi (\gamma_1)}{p}\,(\gamma-1)\,\beta
$$
where 
$$
\beta \,=\,(-1)^{h-1}t^h \partial^h (a\log (\pi))\,=\,(-1)^{h-1}a t^h \partial^{h-1}
\left (\frac{1+\pi}{\pi}\right ).
$$
It implies immediately that $\beta \in F_0\Ddagrig (V).$ On the other hand
$
D^{\Ph=p^{-1}}\,=\,\Cal D_{\text{\rm cris}}(\text{\rm gr}_0 \Ddagrig (V)).
$
Write  $\tilde a$ for  the image of $a$ in $\text{\rm gr}_0\Ddagrig (V)\,\left [1/t \right ]$
and $e_m$ for the canonical base of $D_m.$
Since
$
\text{\rm gr}_0 \Ddagrig (V)\,\simeq \,\underset{i=1}\to{\overset{e}\to \oplus} D_{m_i},
$
without lost of generality  we may assume that $\tilde a\,=\,t^{-m_i}e_{m_i}$ for some $i.$ Let $\tilde \beta$ be
 the image of $\beta$ in $\text{\rm gr}_0 \Ddagrig (V)^{\psi=1}$ and let
$
h^1_0 \,\,:\,\, \text{\rm gr}_0 \Ddagrig (V)^{\psi=1}@>>>H^1(\text{\rm gr}_0\Ddagrig (V))
$
be the canonical map furnished by Proposition 1.1.7.  Recall that $h^1_0(\tilde \beta )=\cl (c,\tilde\beta)$ where 
$(1-\gamma)\,c\,=\,(1-\Ph)\,\tilde \beta.$ 
Then $\tilde \beta= (-1)^{h-1}t^{h-m_i}\partial^h \log (\pi)$. By Lemma 1.5.1 of \cite{CC1}
there exists a unique $b_0\in \boB^{\dag, \psi=0}_{\Bbb Q_p}$ 
such that $(\gamma-1)\,b_0=\ell (X)$. This implies that 
 $$
(1-\gamma)\,(t^{h-m_i}\partial^h b_0e_{m_i})\,=\,(1-\Ph)\,(t^{h-m_i}\partial^h \log (\pi)e_{m_i})\,=\,(-1)^{h-1}
(1-\Ph)\,\tilde \beta.
$$
Thus $c= (-1)^{h-1}t^{h-m_i}\partial^h b_0e_{m_i}$ and 
$
\res (ct^{m_i-1}dt)=(-1)^{h-1}\res (t^{h-1}\partial^h b_0dt)\,e_{m_i}=0.
$
Next from the congruence
$
\tilde \beta \equiv (h-1)!\,t^{-m_i}e_{m_i}\,
\pmod{\Bbb Q_p[[\pi]]\,e_{m_i}}.
$
it follows  that
$\res (\tilde \beta t^{m_i-1}dt)\,=\,(h-1)!\,e_{m_i}.$ Therefore by \cite{Ben2}, Corollary 1.5.6  we have 
$$
\cl (c,\tilde \beta)\,=\,(h-1)!\,\cl (\beta_m)\,=\,(h-1)! \dsize \frac{p}{\log \chi (\gamma_1)}\,i_{\text{\rm gr}_0(\Ddagrig(V)),c}(a).
$$
 On the other hand 
$$
\alpha (0)\,=\,\left.
a \otimes \ell \left (\frac{(1+X)^{\chi (\gamma)}-1}{X}\right )
 \right \vert_{X=0}\,=\,a\left (1-\frac{1}{p} \right )\,\log (\chi (\gamma)).
$$
Theses formulas imply that
$$
\delta_{D,h}(\alpha)=(h-1)!\,\left (1-\frac{1}{p}\right )^{
 -1}(\log \chi
(\gamma))^{-1}\, i_{\text{\rm gr}_0(\Ddagrig(V)),c} (\alpha).
$$
and the proposition is proved.
  \enddemo
\flushpar
{\bf 3.2.3.}  Define
$$
H^1_{f,\{p\}}(V)=\ker \left (H^1_S(V)@>>>\underset{v\in S-\{\infty\}}\to \bigoplus 
\frac{H^1(\Bbb Q_v,V)}{H^1_f(\Bbb Q_v,V)}\right ).
$$

From the definition of $H^1_S(D,V)$ we immediately obtain isomorphisms
$$
\frac{H^1(\Bbb Q_p,V)}{H^1_{f,\{p\}}(V)+ H^1(F_{-1}\Ddagrig (V))}\simeq 
\frac{H^1(F_0\Ddagrig (V))}{H^1_{S}(D,V)+ H^1(F_{-1}\Ddagrig (V))} \simeq \frac{H^1(\text{\rm gr}_0\Ddagrig (V))}{H^1_S(D,V)}\,.
$$ 
Thus, the map $\delta_{D,h}$ constructed in Proposition 3.3.2 induces a map
$$
D^{\Ph=p^{-1}} @>>> 
\dsize\frac{H^1(\Bbb Q_p,V)}{H^1_{f,\{p\}}(V) +H^1(F_{-1}\Ddagrig (V))}
$$
which we will denote again by $\delta_{D,h}.$  On the other hand, we have isomorphisms
$$
D^{\Ph=p^{-1}} \overset{\exp_{V}}\to \iso \frac{H^1_f(\Bbb Q_p,V)}{\exp_{V,\Bbb Q_p}(D_{-1})}
\simeq  \frac{H^1_f(\Bbb Q_p,V)}{H^1(F_{-1}\Ddagrig (V))}
\simeq
\frac{H^1(\Bbb Q_p,V)}{H^1_{f,\{p\}}(V) +H^1(F_{-1}\Ddagrig (V))}.
$$
\proclaim{Proposition 3.2.4} Let $\lambda_D \,:\,D^{\Ph=p^{-1}}@>>> D^{\Ph=p^{-1}}$ denote the homomorphism
making the diagram
$$\xymatrix{
D^{\Ph=p^{-1}} \ar[dr]^{\delta_{D,h}} \ar[rr]^{\lambda_D}
& & D^{\Ph=p^{-1}} \ar[dl]_{(h-1)!\exp_{V}\,\,\,\,}\\
 & {\dsize\frac{H^1(\Bbb Q_p,V)}{H^1_{f,\{p\}}(V)+H^1(F_{-1}\Ddagrig (V))}} &
 }
$$
commute. 
Then 
$$
\det \left (\lambda_D\,\vert \, D^{\Ph=p^{-1}} \right )\,=\,(\log \chi (\gamma))^{-e}\left ( 1-\frac{1}{p} \right )^{-e}
 \Cal L(D,V).
$$
\endproclaim
\demo{Proof}  The proposition follows from Proposition 2.1.6, Proposition 3.2.2 and the following 
elementary fact. Let $U=U_1\oplus U_2$ be the decomposition of a vector space $U$ of dimension
$2e$ into the direct sum of two subspaces of dimension $e$. Let $W\subset U$ be a subspace
of dimension $e$ such that $W\cap U_1=\{0\}.$ Consider the  diagrams 
$$
\xymatrix{
W \ar[r]^{p_1} \ar[d]_{p_2} &U_1  &  U/W  &  U_1 \ar[l]_{i_1} \\
U_2 \ar[ur]_f &  &  U_2 \ar[u]^{i_2} \ar[ur]_g & }
$$
where $p_k$ and $i_k$ are induced by natural projections and inclusions.
Then $f=-g.$ Applying this remark to $U=H^1(\text{\rm gr}_0\Ddagrig (V)),$
$W=H^1_S(D,V)$, $U_1=H^1_f(\text{\rm gr}_0\Ddagrig (V)),$ 
$U_2=H^1_c(\text{\rm gr}_0\Ddagrig (V))$ and taking determinants we obtain the proposition. 
\enddemo

\head {\bf \S4. Special values of $p$-adic $L$-functions}
\endhead

\flushpar
{\bf 4.1. The Bloch-Kato conjecture} (see \cite{F3}, \cite{FP},\cite{BF}).

{\bf 4.1.1.} Let $V$ be a $p$-adic pseudo-geometric representation
of $\text{\rm Gal} (\overline{\Bbb Q}/\Bbb Q).$  Thus $V$ is a
finite-dimensional $\Bbb Q_p$-vector space equipped with a
continuous action of the Galois group $G_S$ for a suitable finite
set of places  $S$ containing $p$. 
Write
$
\RG_S (V)= C_c^\bullet (G_S,V)
$
and define
$$
\RG_{S,c}(V)=\text{\rm cone} \left (\RG_S(V)@>>> \underset{v\in S\cup \{\infty\}}
\to \oplus \RG (\Bbb Q_v,V) \right )\,[-1].
$$

Fix a $\Bbb Z_p$-lattice $T$  of $V$ stable under the action of
$G_S$ and set 
$
\Delta_S(V)={\det}^{-1}_{\Bbb Q_p} \RG_{S,c}(V)
$
and $\Delta_S(T)={\det}^{-1}_{\Bbb Z_p} \RG_{S,c}(T).$
Then $\Delta_S(T)$ is a $\Bbb Z_p$-lattice of the one-dimensional
$\Bbb Q_p$-vector space $\Delta_S(V)$ which does not depend on the
choice of $T$. Therefore it defines a
$p$-adic norm on $\Delta_S(V)$ which we denote by $\Vert \,\cdot
\,\Vert_S.$ Moreother, $(\Delta_S(V), \Vert \,\cdot \,\Vert_S)$
does not depend on the choice of $S.$ More precisely, if $\Sigma$
is a finite set of places which contains $S$, then there exists a
natural isomorphism $\Delta_S(V) @>>>\Delta_{\Sigma}(V)$ such that
$\Vert \,\cdot \,\Vert_\Sigma=\Vert \,\cdot \,\Vert_S.$ It allows
to define the Euler-Poincar\'e line $\Delta_{\text{\rm
EP}}(V)$ as $(\Delta_S(V),\Vert \,\cdot \,\Vert_S)$ where
$S$ is sufficiently large. Recall that for any finite place $v\in
S$ we defined
$$
\RG_f (\Bbb Q_v,V)\,=\,\cases \left [ V^{I_v} @>1-f_v>> V^{I_v}
\right ] &{\text{\rm if $v\ne p$}}\\
\left [\Dc (V)@>(\text{\rm pr},1-\Ph)>> t_V(\Bbb Q_p)\oplus \Dc
(V) \right ] &{\text{\rm if $v=p$}}.
\endcases
$$
At $v=\infty$  we set
$
\RG_f(\Bbb R,V)\,=\,\left [V^+@>>>0\right ],
$
where the first term is placed in degree $0.$ Thus $\RG_f(\Bbb R,V) \overset{\sim}\to \rightarrow \RG(\Bbb R,V)$.
For any $v$ we have a canonical morphism $\text{\rm loc}_p\,:\,\RG_f(\Bbb Q_v,V)@>>> \RG (\Bbb Q_v,V)$
which can
be viewed as a local condition in the sense of \cite{N2}. Consider
the diagram
$$
\xymatrix{ \RG_S(V) \ar[r] & \underset{v\in S\cup \{\infty\}}\to \oplus
\RG(\Bbb
Q_v,V)\\
 & \underset{v\in S\cup \{\infty\}}\to \oplus \RG_f(\Bbb Q_v,V) \ar[u]
 }
$$
and define
$$
\RG_f(V)=\text{\rm cone} \left (\RG_S(V)\oplus \left
(\underset{v\in S\cup\{\infty\}}\to \oplus \RG_f(\Bbb Q_v,V) \right ) @>>>
\underset{v\in S\cup\{\infty\}}\to\oplus \RG (\Bbb Q_v,V)\right )[-1].
$$
Thus, we have a distinguished triangle
$$
\RG_f(V) @>>>\RG_S(V)\oplus \left
(\underset{v\in S\cup\{\infty\}}\to \oplus \RG_f(\Bbb Q_v,V) \right ) @>>>
\underset{v\in S\cup \{\infty\}}\to\oplus \RG (\Bbb Q_v,V).
\tag{4.1}
$$
Set
$$
\Delta_f (V)={\det}_{\Bbb Q_p}^{-1}\RG_f(V)\otimes
{\det}_{\Bbb Q_p}^{-1}t_V(\Bbb Q_p)\otimes {\det}_{\Bbb Q_p}V^+.
$$
\flushpar
It is easy to see that $\RG_f(V)$ and $\Delta_f (V)$ do
not depend on the choice of $S$. Consider the distinguished triangle
$$
\RG_{S,c}(V)@>>>\RG_f(V)@>>>\underset{v\in S\cup \{\infty\}}\to\oplus \RG_f(\Bbb
Q_v,V).
$$
Since ${\det}_{\Bbb Q_p}\RG_f (\Bbb Q_p,V) \simeq {\det}_{\Bbb Q_p}^{-1}t_V(\Bbb Q_p)$ and
${\det}_{\Bbb Q_p}\RG_f(\Bbb R,V)={\det}_{\Bbb Q_p}V^+$
tautologically, we obtain canonical isomorphisms
$$
\Delta_f(V)\simeq {\det}_{\Bbb
Q_p}^{-1}\RG_{S,c}(V)\simeq \Delta_{\text{\rm EP}}(V).
$$
\flushpar
The cohomology of $\RG_f(V)$ is as follow:
$$
\aligned &\bold R^0\Gamma_f(V)=H_S^0(V),\quad
\bold R^1\Gamma_f(V)=H^1_f(V),\quad
\bold R^2\Gamma_f(V)\simeq H^1_f(V^*(1))^*,\\
&\bold R^3\Gamma_f (V)=\text{\rm coker}\,\left (
H^2_S(V)@>>>\underset{v\in S}\to\oplus H^2(\Bbb Q_v,V)\right )
\simeq H_S^0(V^*(1))^*.
\endaligned
\tag{4.2}
$$
These groups seat in the following exact sequence:
$$
\multline 0@>>>\bold R^1\Gamma(V)@>>>H^1_S(V)@>>>\underset{v\in
S}\to \bigoplus \frac{H^1(\Bbb Q_v,V)}{H^1_f(\Bbb Q_v,V)}
@>>>\bold
R^2\Gamma_f(V)@>>>\\
H^2_S(V)@>>>\underset{v\in S}\to \oplus H^2(\Bbb Q_v,V)@>>>\bold
R^3\Gamma_f(V)@>>>0.
\endmultline
$$
\flushpar
The $L$-function of $V$ is defined as the Euler product
$$
L(V,s)=\underset{v}\to \prod E_v(V,(Nv)^{-s})^{-1}
$$
where
$$
E_v(V,t)=\cases \det \left (1-f_vt\,\vert \,V^{I_v} \right ),
&\text{\rm if $v\ne p$} \\
\det \left (1-\Ph t\,\vert \,\Dc (V)\right ) &\text{\rm if $v=p$}.
\endcases
$$

$\,$
\newline
\flushpar
{\bf 4.1.2.} In this paper we treat motives in
the formal sense and assume all conjectures about the category of
mixed motives $\Cal M \Cal M$ over $\Bbb Q$ which are  necessary to state the Bloch-Kato
conjecture (see \cite{F3}, \cite{FP}). If $M$ is a pure motive
over $\Bbb Q$ we denote by $M_v$ its $v$-adic realizations. Assume that the groups
$
H^i(M)\,=\,\text{\rm Ext}^i_{\Cal M\Cal M}(\Bbb Q(0),M)
$
are well defined and vanish for $i\ne 0,1.$ It should be possible to  define  a
$\Bbb Q$-subspace $H^1_f(M)$ of $H^1(M)$ consisting of "integral"
classes of extensions which is expected to be finite dimensional.
It is convenient to set $H^0_f(M)=H^0(M).$ Then we assume that for
any finite place $v$  the regulator map induces   isomorphisms
$$
H^i_f(M)\otimes_{\Bbb Q}\Bbb Q_p\simeq H^i_f(M_v), \qquad i=0,1. \tag{4.3}
$$

Let $M$ be a motive satisfying the following condition
\newline
\,

{\bf M)} $H^i_f(M)=H^i_f(M^*(1))$ for $i=0,1.$
\newline
\,
\flushpar
Let $M_{\text{\rm dR}}$ and $M_{\text{\rm B}}$ denote the de Rham
and the Betti realizations of $M$ respectively and let $t_M(\Bbb
Q)= M_{\text{\rm dR}}/\F^0 M_{\text{\rm dR}}$ be the tangent space
of $M.$ The complex conjugation $c$ acts on $M_{\text{\rm B}}$ and 
$
M_{\text{\rm B}}=M_{\text{\rm B}}^+\oplus M_{\text{\rm B}}^-.
$
The  comparision isomorphism
$
M_{\text{\rm B}}\otimes_{\Bbb Q}\Bbb R \simeq M_{\text{\rm
dR}}\otimes_{\Bbb Q} \Bbb R
$
induces a map
$$
M_{\text{\rm B}}^+\otimes_{\Bbb Q}\Bbb R@>>>t_M(\Bbb R)
$$
which is expected to be an isomorphism. Assuming this, we can
define a natural injective map
$$
\Omega_{\infty}\,:\,{\det}_{\Bbb Q}^{-1}t_M(\Bbb Q) \otimes
{\det}_{\Bbb Q}M_{\text{\rm B}}^+@>>> \Bbb R.
$$
Fix $\omega_t\in {\det}_{\Bbb Q}t_M(\Bbb Q)$ and $\omega_{\text{\rm
B}}\in {\det}_{\Bbb Q} M_{\text{\rm B}}^+$ and set
$
\Omega_{\infty}(\omega_t, \omega_{\text{\rm B}})=\Omega_{\infty}(\omega_t^{-1}\otimes \omega_{\text{\rm
B}}).
$
\flushpar
It is conjectured that the $L$-function $L(M_v,s)$ does not depend
of $v$. It will be denoted by $L(M,s).$

\proclaim\nofrills{Conjecture } $\,\,${\smc (Deligne)}. Let $M$ be a
motive satisfying {\bf M)}. Then
$$
\frac{L(M,0)}{\Omega_{\infty}(\omega_t, \omega_{\text{\rm B}})}\in \Bbb Q^*.
$$
\endproclaim

\flushpar
{\bf 4.1.3.} Let $p$ be a prime number and let $M_p$ denote the $p$-adic
realization of $M.$ From {\bf M)} and (4.3) it follows that
$H^0(M_p)=H^0(M_p^*(1))=0$ and $H^1_f(M_p)=H^1_f(M_p(1))=0.$ Hence
$\RG_f (M_p)$ is acyclic. Fix $\omega_t$ and $\omega_{\text{\rm
B}}$ and define a  map
$$
i_{\omega_t,\omega_{\text{\rm
B}},p}\,\,:\,\,\Delta_{\text{\rm EP}}(M_p)\overset{\sim}\to \rightarrow
{\det}_{\Bbb Q_p}^{-1} t_M(\Bbb Q_p) \otimes {\det}_{\Bbb Q_p}
M_{\text{\rm B}}^+ @>>> \Bbb Q_p
$$
by
$
x=i_{\omega_t,\omega_{\text{\rm
B}},p} (x)\,(\omega_t^{-1}\otimes \omega_{\text{\rm B}}).
$
The Bloch-Kato conjecture states as follow:

\proclaim\nofrills{Conjecture} $\,\,${\smc (Bloch-Kato)}. Let $T_p$ be a
$\Bbb Z_p$-lattice of $M_p$ stable under the action of $G_S.$ Then
$$
i_{\omega_t,\omega_B,p}(\Delta_{\text{\rm
EP}}(T_p))\,=\,\frac{L(M,0)}{\Omega_{\infty}(\omega_t, \omega_{\text{\rm B}})}\,\Bbb
Z_p.
$$
\endproclaim

\flushpar
{\bf 4.2. The complex $\bold R\Gamma^{(\eta_0)}_{\text{\rm Iw},h} (D,V)$.}
\newline
{\bf 4.2.1.} 
Let $\Gamma$ denote the Galois group of $\Bbb Q(\zeta_{p^\infty})/\Bbb Q$ and 
$\Gamma_n=\Gal (\Bbb Q(\zeta_{p^\infty})/\Bbb Q (\zeta_{p^n})).$ Set
$\Lambda=\Bbb Z_p[[\Gamma_1]]$ and 
$\Lambda (\Gamma)= \Bbb Z_p[\Delta]\otimes_{\Bbb Z_p} \Lambda$.
For any character $\eta \in X(\Delta)$ put
$$
e_{\eta}\,=\,\frac{1}{|\Delta|}\sum_{g\in \Delta} \eta^{-1}(g)g.
$$
Then $\Lambda (\Gamma) =\underset{\eta \in X(\Delta)}\to \oplus \Lambda (\Gamma)^{(\eta)}$ where
$\Lambda (\Gamma)^{(\eta)}=\Lambda e_\eta$ and for any $\Lambda (\Gamma)$-module 
$M$ one has a canonical decomposition 
$$
M\simeq \oplus_{\eta \in X(\Delta)}M^{(\eta)},\qquad M^{(\eta)}=e_\eta (M).
$$ 
We write $\eta_0$ for the trivial character of $\Delta$ and 
identify $\Lambda$ with $\Lambda (\Gamma)e_{\eta_0}.$
\newline
\newline
Let $V$ be a p-adic pseudo-geometric representation  
unramified outside $S.$ Set $d(V)=\dim (V)$ and $d_{\pm}(V)=\dim (V^{c=\pm 1}).$
Fix a $\Bbb Z_p$-lattice $T$ of $V$ stable under the action of $G_{S}.$
Let $\iota \,:\,\Lambda (\Gamma)@>>>\Lambda (\Gamma)$ denote the canonical involution $g\mapsto g^{-1}.$ Recall that
the induced module $\text{\rm Ind}_{\Bbb Q(\zeta_{p^\infty})/\Bbb Q} (T)$
is isomorphic to $(\Lambda (\Gamma)\otimes_{\Bbb Z_p}T)^\iota $ (\cite{N2}, section 8.1).
Define
$$
\align
&H^i_{\text{\rm Iw},S}(T)\,=\, 
H^i_S((\Lambda (\Gamma)\otimes_{\Bbb Z_p}T)^\iota),\\
&\Hi^i(\Bbb Q_v,T)\,=\,H^i(\Bbb Q_v,(\Lambda (\Gamma)\otimes_{\Bbb Z_p}T)^\iota)
\qquad \text{\rm for any finite place $v$.}
\endalign
$$
From Shapiro's lemma it follows immediately that
$$
H^i_{\text{\rm Iw},S}(T)\,=\,\varprojlim_{\text{\rm cores}} 
H^i_S(\Bbb Q(\zeta_{p^n}),T),\qquad
\Hi^i(\Bbb Q_p,T)\,=\,\varprojlim_{\text{\rm cores}} H^i(\Bbb Q_p(\zeta_{p^n}),T).
$$
Set $H^i_{\text{\rm Iw},S}(V)\,=\,H^i_{\text{\rm Iw},S}(T)\otimes_{\Bbb Z_p}\Bbb Q_p$
and $H^i_{\text{\rm Iw}}(\Bbb Q_v,V)\,=\,
H^i_{\text{\rm Iw}}(\Bbb Q_v,T)\otimes_{\Bbb Z_p}\Bbb Q_p.$ In \cite{PR2}
Perrin-Riou   proved the following results about the structure of these
modules. 
\newline
\,

i)  $H^i_{\text{\rm Iw},S}(V)=0$ and $H^i_{\text{\rm Iw}}(\Bbb Q_v,T)=0$ if $i\ne 1,2;$
\newline
\,

ii) If $v\ne p$, then for each $\eta\in X(\Delta)$ the $\eta$-component
 $H^i_{\text{\rm Iw}}(\Bbb Q_v,T)^{(\eta)}$ is a finitely generated
torsion $\Lambda$-module. In particular, 
 $H^1_{\text{\rm Iw}}(\Bbb Q_v,T)\simeq 
H^1(\Bbb Q_v^{\text{ur}}/\Bbb Q_v,(\Lambda (\Gamma)\otimes_{\Bbb Z_p}T^{I_v})^\iota)$.
\newline
\,

iii) If $v=p$ then $H^2_{\text{\rm Iw}}(\Bbb Q_p,T)^{(\eta)}$ are finitely generated torsion $\Lambda$-modules. 
Moreover, for each $\eta\in X(\Delta)$

$$
\text{\rm {rg}}_{\Lambda} \left (H^1_{\text{\rm Iw}}(\Bbb Q_p,T)^{(\eta)}\right )\,=\,d,\qquad
H^1_{\text{\rm Iw}}(\Bbb Q_p,T)^{(\eta)}_{\text{\rm tor}}\simeq H^0(\Bbb Q_p(\zeta_{p^\infty})\,,T)^{(\eta)}.
$$
Remark that by local duality $H^2_{\text{\rm Iw}}(\Bbb Q_p,T)\simeq H^0(\Bbb Q_p(\zeta_{p^{\infty}}),V^*(1)/T^*(1))$.
\newline
\,

iv) If  the weak Leopoldt conjecture holds for the pair $(V,\eta)$  i.e. if 
$H^2_S(\Bbb Q(\zeta_{p^\infty}),V/T)^{(\eta)}=0$
then $H^2_{\text{\rm Iw},S}(T)^{(\eta)}$ is $\Lambda$-torsion and  
$$
\text{\rm rank}_{\Lambda}\left ( H^1_{\text{\rm Iw},S}(T)^{(\eta)}\right )\,=\,
\cases d_-(V), &{\text{\rm if $\eta (c)=1$}}\\
d_+(V), &{\text{\rm if $\eta (c)=-1$.}}
\endcases
$$ 

Passing to the projective limit in the Poitou-Tate exact sequence one obtains
an exact sequence
$$
\multline
0@>>> H^2_S(\Bbb Q(\zeta_{p^\infty}), V^*(1)/T^*(1))^{\wedge}@>>>
H^1_{\text{\rm Iw},S}(T)@>>> \underset{v\in S}\to \oplus H^1_{\text{\rm Iw}}(\Bbb Q_v,T)@>>>
H^1_S(\Bbb Q(\zeta_{p^\infty}), V^*(1)/T^*(1))^{\wedge}\\
@>>>H^2_{\text{\rm Iw},S}(T)@>>> \underset{v\in S}\to \oplus H^2_{\text{\rm Iw}}(\Bbb Q_v,T)@>>>
H^0_S(\Bbb Q(\zeta_{p^\infty}), V^*(1)/T^*(1))^{\wedge}@>>>
0.
\endmultline
\tag{4.4}
$$
Define
$$
\align
&\RG_{\text{\rm Iw},S} (T)\,=\,C_c^\bullet (G_S,(\Lambda (\Gamma)\otimes_{\Bbb Z_p} T)^\iota),\\
&\RG_{\text{\rm Iw}}(\Bbb Q_v,T)\,=\,C_c^\bullet (G_v,(\Lambda (\Gamma)\otimes_{\Bbb Z_p} T)^\iota),\\
&\RG_S(\Bbb Q(\zeta_{p^\infty}), V^*(1)/T^*(1))\,=\,
C_c^\bullet (G_{S},\text{\rm Hom}_{\Bbb Z_p}(\Lambda (\Gamma), V^*(1)/T^*(1))).
\endalign
$$
Then 
the sequence (4.3) is induced by the distinguished triangle
$$
\RG_{\text{\rm Iw},S} (T)@>>>\underset{v\in S}\to \oplus \RG_{\text{\rm Iw}}(\Bbb Q_v,T)@>>>
\left (\RG_S(\Bbb Q(\zeta_{p^\infty}), V^*(1)/T^*(1))^\iota \right )^{\wedge}\,[-2]
$$
(\cite{N2}, Theorem 8.5.6). Finally, we have  usual descent formulas
$$
\RG_{\text{\rm Iw},S} (T)\otimes_{\Lambda}^{\bold L}\Bbb Z_p \simeq \RG_{S} (T),\qquad 
\RG_{\text{\rm Iw}} (\Bbb Q_v,T)\otimes_{\Lambda}^{\bold L}\Bbb Z_p \simeq \RG (\Bbb Q_v,T)
$$
( \cite{N2}, Proposition 8.4.21).
\newline
\newline
{\bf 4.2.2.} For the remainder of this chapter we assume that $V$ satisfies  the conditions {\bf C1-5)} of section 3.1.2
where {\bf C2)} is replaced  by 
the following  stronger condition
\newline
\,

{\bf C2*)} $H^0(\Bbb Q_p,V)= H^0(\Bbb Q_p, V^*(1))= 0$.
\newline
\,
\flushpar
Remark that {\bf C1)} and {\bf C2*)} guarantee that the weak Leopoldt conjecture holds for $(V,\eta_0)$ and $(V^*(1),\eta_0)$ 
( Proposition B.5 of \cite{PR2}). To simplify notations we write  $\Cal H$ for $\Cal H(\Gamma_1)$.
In this subsection we interpret Perrin-Riou's construction 
of the module of $p$-adic $L$-functions in terms of \cite{N2}.
Fix an admissible subspace $D$ of $\Dc (V)$ and  a $\Bbb Z_p$-lattice $N$ of $D$.
Set $\Cal D_p(N,T)^{(\eta_0)}\,=\,N\otimes_{\Bbb Z_p}\Lambda$, 
$\RG_{\text{\rm Iw},f}^{(\eta_0)}(\Bbb Q_p,N,T) =\Cal D_p(N,T)^{(\eta_0)}[-1]$ and
$\RG_{\text{\rm Iw},f}^{(\eta_0)}(\Bbb Q_p,D,V)=\RG_{\text{\rm Iw},f}^{(\eta_0)}(\Bbb Q_p,N,T)\otimes_{\Bbb Z_p}
\Bbb Q_p.$
Consider the map
$$
\bExp_{V,h}^\ep \,:\,\RG_{\text{\rm Iw},f}^{(\eta_0)}(\Bbb Q_p,T)\otimes_{\La}\Cal H
@>>> \RG_{\text{\rm Iw}}^{(\eta_0)} (\Bbb Q_p,T)\otimes_{\Lambda}^{\bold L}\Cal H
$$
which will be viewed as a local condition at $p$.  If  $v\ne p$ the inertia group $I_v$ acts trivially
on $\Lambda$  set
$$
\RG_{\text{\rm Iw},f}^{(\eta_0)}(\Bbb Q_v,N,T)\,=\,
\left [T^{I_v}\otimes \Lambda^{\iota}  @>1-f_v>>T^{I_v}\otimes \Lambda^{\iota}  \right ] 
$$
where the first term is placed in degree $0$. We have a commutative diagram
$$
\xymatrix{
\RG_{\text{\rm Iw},S}^{(\eta_0)}(T)\otimes_{\Lambda}\Cal H
\ar[r] & \underset{v\in S} \to \oplus \RG_{\text{\rm Iw}}^{(\eta_0)}(\Bbb Q_v,T)\otimes_{\Lambda} \Cal H\\
 & \underset{v\in S}\to \oplus \RG_{\text{\rm Iw},f}^{(\eta_0)}(\Bbb Q_v,N,T)\otimes_{\Lambda} \Cal H
 \ar[u]
}
\tag{4.5}
$$
Consider the associated Selmer complex
$$
\multline
\RG_{\text{\rm Iw},h}^{(\eta_0)}(D,V)\,=\\
\text{\rm cone}\, \left [
\left ( \RG_{\text{\rm Iw},S}^{(\eta_0)}(T)\oplus 
\left (\underset{v\in S}\to\oplus \RG_{\text{\rm Iw},f}^{(\eta_0)}(\Bbb Q_v,N,T)\right ) \right )\otimes_{\Lambda} \Cal H
@>>>
\underset{v\in S}\to \oplus \RG_{\text{\rm Iw}}^{(\eta_0)}(\Bbb Q_v,T)\otimes_{\Lambda} \Cal H\right ]
[-1]
\endmultline
$$
It is easy to see that it does not depend on the choice of $S$. Our main result 
about this complex is the following theorem.

\proclaim{Theorem 4.2.3} Assume that $V$ satisfies the conditions C1-5).
Let $D$ be an admissible subspace of $\Dc (V).$ Assume that $\Cal L(V,D)\ne 0.$  
Then

i) $\bold R^i\Gamma_{\text{\rm Iw},h}^{(\eta_0)}(D,V)$ are $\Cal H$-torsion modules for all $i$.

ii)  $\bold R^i\Gamma_{\text{\rm Iw},h}^{(\eta_0)}(D,V)=0$ for  $i\ne 2,3$ and
$$
\bold R^3\Gamma_{\text{\rm Iw},h}^{(\eta_0)}(D,V)\simeq \left(H^0(\Bbb Q(\zeta_{p^{\infty}}),V^*(1))^*
\right )^{(\eta_0)}\otimes_{\Lambda}\Cal H.
$$

iii) The complex $ \bold R\Gamma_{\text{\rm Iw},h}^{(\eta_0)}(D,V)$ is semisimple i.e.
for each $i$ the natural map 
$$
\bold R^i\Gamma_{\text{\rm Iw},h}^{(\eta_0)}(D,V)^{\Gamma} @>>>
\bold R^i\Gamma_{\text{\rm Iw},h}^{(\eta_0)}(D,V)_{\Gamma}
$$
is an isomorphism.
\endproclaim
\flushpar
{\bf 4.2.4. Proof of Proposition 4.2.3.} We leave the proof
of the following lemma as an easy exercise. 

\proclaim{Lemma 4.2.4.1} Let $A$ and $B$ be two submodules of a finitely generated 
free $\Cal H$-module $M$.
Assume that the  natural maps $A_{\Gamma_1} @>>> M_{\Gamma_1}$ and $B_{\Gamma_1} @>>> M_{\Gamma_1}$ are both
injective. Then $A_{\Gamma_1}\cap B_{\Gamma_1}=\{0\}$ implies that $A\cap B=\{0\}.$
\endproclaim
\flushpar
{\bf 4.2.4.2.} Since $H_{\text{\rm Iw},S}^0(V)$ and
$\Hi^0(\Bbb Q_v,V)$ are zero, we have   $ \bold R^0\Gamma_{\text{\rm Iw},h}^{(\eta_0)}(D,V)=0.$
Next, by definition $\bold R^1\Gamma_{\text{\rm Iw},h}(D,V)^{(\eta_0)}\,=\,\ker (f)$ where
$$
f\,:\, \left (
H^1_{\text{\rm Iw},S}(T)^{(\eta_0)}\oplus 
\Cal D_p(N,T)^{(\eta_0)}
\underset{v\in S-\{p\}}\to \oplus H^1_{\text{\rm Iw},f}(\Bbb Q_v,T)^{(\eta_0)}
\right )\otimes \Cal H @>>>
\underset{v\in S}\to \oplus \Hi^1(\Bbb Q_v,T)^{(\eta_0)}\otimes \Cal H 
$$ 
is the map induced by (4.5).
If $v\in S-\{p\}$ one has   
$$
H^1_{\text{\rm Iw},f}(\Bbb Q_v,T)^{(\eta_0)}\,=\,\Hi^1(\Bbb Q_v,T)^{(\eta_0)}= 
H^1(\Bbb Q_v^{\text{\rm ur}}/\Bbb Q_v, (\Lambda\otimes T^{I_v})^{\iota}).
$$  
Thus
$$
\bold R^1\Gamma_{\text{\rm Iw},h}^{(\eta_0)}(D,V)\,=\,
\left (H^1_{\text{\rm Iw},S}(T)^{(\eta_0)}\otimes_{\Lambda}\Cal H \right )  \cap \left (\Exp^\ep_{V,h}
\left (\Cal D_p(D,T)^{(\eta_0)}\right  )\otimes_{\Lambda} \Cal H\right ) 
$$ 
in  $\Hi^1(\Bbb Q_p,T)^{(\eta_0)}\otimes_{\Lambda}\Cal H$. Put
$$
A=\Exp_{V,h}^\ep (D_{-1}\otimes \Cal H)\oplus X^{-1} \Exp_{V,h}^\ep (D^{\Ph=p^{-1}}\otimes \Cal H)
\subset \Hi^1(\Bbb Q_p,T)^{(\eta_0)}\otimes_{\Lambda}\Cal H.
$$
By Theorem 2.2.4 and Proposition 3.2.2 $A_{\Gamma_1}$ injects into 
$H^1(\Bbb Q_p,V).$ The $\Cal H$-module  $M=\left (\dsize \frac{\Hi^1(\Bbb Q_p,T)}{T^{H_{\Bbb Q_p}}}\right )^{(\eta_0)}\otimes_{\Lambda}\Cal H $ is  free and $A\hookrightarrow M.$ Since  
$T^{G_{\Bbb Q_p}}=0$ by {\bf C2*)}, one has
$
M_{\Gamma_1}\,=\,\Hi^1(\Bbb Q_p,V)_{\Gamma}\subset  H^1(\Bbb Q_p,V)
$
and we obtain that  $A_{\Gamma_1}$ injects into $M_{\Gamma_1}.$

Set $B=\left (\dsize \frac{H^1_{\text{\rm Iw},S}(T)}{T^{H_{\Bbb Q}}} \right )^{(\eta_0)}\otimes_{\Lambda}\Cal H.$ 
The weak Leopoldt conjecture for $(V^*(1),\eta_0)$ together with
the fact that $\Hi^1(\Bbb Q_v,T)$ are $\Lambda$-torsion for $v\in S-\{p\}$
imply that  $B\hookrightarrow M.$  Since the image of $\Hi^1(\Bbb Q_v,V)_{\Gamma}$ in
$H^1(\Bbb Q_v,V)$ is contained in $H^1_f(\Bbb Q_v,V),$ the image of $H^1_{\text{\rm Iw},S}(V)_{\Gamma}$
in $H^1_S(V)$ is in fact contained in
$$
H^1_{f,\{p\}}(V)\,=\,\ker \left (H^1_S(V)@>>>\underset{v\in S-\{p\}}\to \bigoplus \frac{H^1(\Bbb Q_v,V)}{H^1_f(\Bbb Q_v,V)}
\right ).
$$
Because $H^1_f(V)=0,$ the group $H^1_{f,\{p\}}(V)$ injects into $H^1(\Bbb Q_p,V)$ and we have 
$$
H^1_{\text{\rm Iw},S}(V)^{(\eta_0)}_{\Gamma_1}=H^1_{\text{\rm Iw},S}(V)_{\Gamma}\hookrightarrow H^1_{f,\{p\}}(V)\hookrightarrow
H^1(\Bbb Q_p,V).
$$ 
Thus $B_{\Gamma_1}\subset M_{\Gamma_1}.$ We shall prove that $\bold R^1\Gamma_{\text{\rm Iw},h^{(\eta_0)}}(D,V)=0.$
By Lemma 4.2.4.1 it suffices to show that $A_{\Gamma_1} \cap B_{\Gamma_1}=\{0\}.$
Now we  claim that $ A_{\Gamma_1}\cap H^1_{f,\{p\}}(V)=\{0\}.$ 
First remark that
$$
 H^1_{f,\{p\}}(V)\hookrightarrow \frac{H^1(\Bbb Q_p,V)}{H^1(F_{-1}\Ddagrig (V))}.
$$ 
On the other hand, from Theorem 2.2.4 
it follows that 
$$
\Exp_{V,h}^\ep(D_{-1}\otimes \Cal H)_{\Gamma_1}=\exp_{V,\Bbb Q_p} (D_{-1})\subset H^1(F_{-1}\Ddagrig (V)).
$$
Therefore, Proposition 3.2.2 implies that the image of $A_{\Gamma_1}$ in   
$\dsize \frac{H^1(\Bbb Q_p,V)}{H^1(F_{-1}\Ddagrig (V))}$  coincides with
$H^1_c(\text{\rm gr}_0 \Ddagrig (V)).$ But $\Cal L(D,V)\ne 0$
if and only if 
$H^1_S(D,V)\cap H^1_c(\text{\rm gr}_0 \Ddagrig (V)) =0
$ 
where
$H^1_S(D,V)$ denotes  the inverse image of $H^1(\text{\rm gr}_0 \Ddagrig (V))$ in $H^1_{f,\{p\}}(V).$
This proves the claim and implies that $\bold R^1\Gamma_{\text{\rm Iw},h}^{(\eta_0)}(D,V)=0.$

{\bf 4.2.4.3.} We shall show that $\bold R^2\Gamma_{\text{\rm Iw},h}^{(\eta_0)}(D,V)$ is $\Cal H$-torsion. 
By definition, we have an exact sequence  
$$
0@>>> \text{\rm coker} (f)@>>>\bold R^2\Gamma_{\text{\rm Iw},h}^{(\eta_0)}(D,V)@>>>\sha^2_{\,\text{\rm Iw},S}(V)^{(\eta_0)}\otimes_{\Lambda_{\Bbb Q_p}}
\Cal H@>>>0,
\tag{4.6}
$$ 
where
$$
\sha^2_{\,\text{\rm Iw},S}(V)\,=\,\ker \left ( H^2_{\text{\rm Iw},S}(V)
@>>>
\underset{v\in S}\to \oplus \Hi^2(\Bbb Q_v,V) \right ).
$$
It follows from the weak Leopoldt conjecture that $\sha^2_{\,\text{\rm Iw},S}(V)$ is $\Lambda_{\Bbb Q_p}$-torsion.
On the other hand, as $\Cal H$ is a Bezout ring \cite{La}, the formulas
$$
\text{\rm rank}_{\Lambda} H^1_{\text{\rm Iw},S} (T)^{(\eta_0)}=d_{-}(V),\quad
\text{\rm rank}_{\Lambda} H^1_{\text{\rm Iw}} (\Bbb Q_p,T)^{(\eta_0)}=d(V),\quad
\text{\rm rank}_{\Lambda} \Cal D_p(N,T)=d_{+}(V)
$$
together with the fact that $\bold R^1\Gamma_{\text{\rm Iw},h}^{(\eta_0)}(D,V)=0$ 
imply that $\text{\rm coker} (f)$ is $\Cal H$-torsion. We have therefore proved that
$\bold R^2\Gamma_{\text{\rm Iw},h}(D,V)$ is $\Cal H$-torsion. Finally,
 the Poitou-Tate exact sequence  gives that 
$$
\bold R^3\Gamma_{\text{\rm Iw},h}^{(\eta_0)}(D,V)= \left (H^0(\Bbb Q(\zeta_{p^\infty}),V^*(1))^*\right )^{(\eta_0)}
\otimes_{\Lambda_{\Bbb Q_p}}\Cal H
$$
is also $\Cal H$-torsion. The proposition is proved.
\newline
\,

{\bf 4.2.4.4.} Now we prove the semisimplicity of 
$\bold R \Gamma_{\text{\rm Iw},h}^{(\eta_0)}(D,V).$
 First, remark that {\bf C2*)} implies that  $\Hi^1(\Bbb Q_p,V)^{\Gamma}=0$ and  
$\Hi^1(\Bbb Q_p,V)_{\Gamma}=H^1(\Bbb Q_p,V).$
Next, $H^1_{\text{\rm Iw},S}(V)^{(\eta_0)}\simeq \Lambda_{\Bbb Q_p}^{d_{-}(V)}\oplus 
H^1_{\text{\rm Iw},S}(V)^{(\eta_0)}_{\text{\rm tor}}. $ 
Since  $H^1_{\text{\rm Iw},S} (V)_{\text{\rm tor}}
\subset V^{H_{\Bbb Q_p}},$ we have  
$(H^1_{\text{\rm Iw},S} (V)_{\text{\rm tor}})_{\Gamma}=0$ by the snake
lemma. Thus $\dim_{\Bbb Q_p} H^1_{\text{\rm Iw},S} (V)^{(\eta_0)}_{\Gamma_1}=d_{-}(V).$ 
On the other hand $\dim_{\Bbb Q_p} H^1_{f,\{p\}}(V)=\dim_{\Bbb Q_p} H^1(\Bbb Q_p,V)-
\dim_{\Bbb Q_p} t_V=d_{-}(V).$ Since $H^1_{\text{\rm Iw},S} (V)^{(\eta_0)}_{\Gamma_1}$ injects
into $H^1_{f,\{p\}}(V)$ this proves that $H^1_{\text{\rm Iw},S} (V)^{(\eta_0)}_{\Gamma_1}=
H^1_{f,\{p\}}(V).$
Consider the exact sequence
$$
0@>>>
\left (
H^1_{\text{\rm Iw},S}(T)^{(\eta_0)}\oplus 
\Cal D_p(N,T)^{(\eta_0)}
\right )\otimes \Cal H @>>>
\Hi^1(\Bbb Q_p,T)^{(\eta_0)}\otimes \Cal H @>>>\text{\rm coker} (f)@>>>0.
$$ 
Recall that $\Exp_{V,h,0}^\ep\,:\,D @>>>\Hi^1 (\Bbb Q_p,V)_{\Gamma}$ denotes the homomorphism
induced by the large exponential map. 
Applying the snake lemma, and 
taking into account that $\text{\rm Im} (\Exp_{V,h,0}^\ep)=
\exp_{V,\Bbb Q_p} (D_{-1})=H^1(F_{-1}\Ddagrig (V))$ and  
$\ker (\Exp_{V,h,0}^\ep)= D^{\Ph=p^{-1}}$ (see for example \cite{BB}, Propositions 4.17 and 4.18
or the proof of Proposition 3.3.2) we obtain
$$
\align &\text{\rm coker} (f)^{\Gamma_1}=\ker \left ( 
H^1_{f,\{p\}} (V)\oplus D@>\Exp^\ep_{V,h,0}>> H^1(\Bbb Q_p,V) \right )=D^{\Ph=p^{-1}},\\
&\text{\rm coker} (f)_{\Gamma_1}=\frac{H^1(\Bbb Q_p,V)}{H^1_{f,\{p\}}(V)+H^1(F_{-1}\Ddagrig (V))}.
\endalign
$$
Thus on has a commutative diagram
$$
\CD
\text{\rm coker} (f)^{\Gamma_1} @>>> D^{\Ph=p^{-1}}\\
@VVV   @VV{\delta_{D,h}}V \\
\text{\rm coker} (f)_{\Gamma_1}  @>>> \dsize\frac{H^1(\Bbb Q_p,V)}{H^1_{f,\{p\}}(V)+H^1(F_{-1}\Ddagrig (V))}.
\endCD
$$
where horizontal arrows are isomorphisms and the left vertical arrow is the natural projection.
From Proposition 3.2.4 it follows that  $\text{\rm coker} (f)^{\Gamma_1}@>>>
\text{\rm coker} (f)_{\Gamma_1}$ is an isomorphism if and only if $\Cal L(D,V)\ne 0.$
\newline
On the other hand, the arguments \cite{PR2}, section 3.3.4 show that 
$\sha^2_{\,\text{\rm Iw},S}(V)_{\Gamma}=\sha^2_{\,\text{\rm Iw},S}(V)^{\Gamma}=0.$
Remark that Perrin-Riou assumes that $\Dc (V)^{\Ph=1}=\Dc (V)^{\Ph=p^{-1}}=0$, but
her proof  works in our case without modifications and we repeat it for the
commodity of the reader. Consider the following commutative diagram
$$
\xymatrix{
\underset{v\in S}\to \oplus \Hi^1(\Bbb Q_v,V)_{\Gamma} \ar[r] \ar[d] 
&\left (H^1_S(\Bbb Q(\zeta_{p^\infty}),V^*(1))^*\right )_{\Gamma} \ar[r] \ar[d]
&\sha^2_{\,\text{\rm Iw},S}(V)_{\Gamma} \ar[r] &0\\
\underset{v\in S-\{p\}}\to \oplus H^1_f (\Bbb Q_v,V)\oplus H^1(\Bbb Q_p,V) \ar[r] 
&H^1_S(V^*(1))^* \ar[r] & 0 &}
$$
The top row of this diagram is obtained by taking coinvariants in the Poitou-Tate
exact sequence. Thus it is  exact. The bottom row is obtained from  the exact sequence
$$
0@>>> H^1_S(V^*(1))@>>>H^1(\Bbb Q_p,V^*(1))\oplus \underset{v\in S-\{p\}}\to \bigoplus
\dsize\frac{H^1(\Bbb Q_v,V^*(1))}{H^1_f(\Bbb Q_v,V^*(1))}
$$
by taking duals. Thus, it is an exact sequence too. 
Since $\Hi^1(\Bbb Q_p,V)_{\Gamma}=H^1(\Bbb Q_p,V)$ and 
$ \Hi^1(\Bbb Q_v,V)_{\Gamma}=H^1_f(\Bbb Q_v,V)$ the left vertical map is an isomorphism. 
The right vertical map seats in the exact sequence
$$
0@>>> \left(H^1_S(\Bbb Q(\zeta_{p^\infty}),V^*(1))^*\right )_{\Gamma}@>>>H^1_S(V^*(1))^*
@>>>\left (H^0_S( \Bbb Q(\zeta_{p^\infty}),V)^*\right )^{\Gamma}@>>>0
$$
(see \cite{PR2}, formula (1.4)). The isomorphism $H^0( \Bbb Q_p(\zeta_{p^\infty}),V)
\simeq \underset{k\in \Bbb Z}\to \oplus V(-k)^{G_{\Bbb Q_p}}\,(k)$ together
with the fact that $V^{G_{\Bbb Q_p}}=0$ implies that 
$H^0_S( \Bbb Q(\zeta_{p^\infty}),V)=0$ and the right vertical arrow of the diagram 
is an isomorphism too. 
This proves  that $\sha^2_{\,\text{\rm Iw},S}(V)_{\Gamma}=0.$
Finally, from  $\text{\rm dim}_{\Bbb Q_p} \sha^2_{\,\text{\rm Iw},S}(V)^{\Gamma}
\leqslant \text{\rm dim}_{\Bbb Q_p} \sha^2_{\,\text{\rm Iw},S}(V)_{\Gamma}$ 
it follows that $\sha^2_{\,\text{\rm Iw},S}(V)^{\Gamma}=0.$
Therefore, applying the snake lemma to (4.6) we obtain a commutative diagram
$$
\xymatrix{
\text{\rm coker} (f)^\Gamma \ar[r] \ar[d] &\bold R^2\Gamma^{(\eta_0)}_{\text{\rm Iw},h}(D,V)^{\Gamma}\ar[d]\\
\text{\rm coker} (f)_\Gamma \ar[r] &\bold R^2\Gamma^{(\eta_0)}_{\text{\rm Iw},h}(D,V)_{\Gamma},}
$$
in which  the horizontal arrows are isomorphisms and the vertical arrows are natural projections.
This proves that $\bold R\Gamma^{(\eta_0)}_{\text{\rm Iw},h}(D,V)$ is semisimple
in degree $2$. Remark that the semisimplicity in degree $3$  is obvious because by ii)
$
\bold R^3\Gamma^{(\eta_0)}_{\text{\rm Iw},h}(D,V)^{\Gamma}=
\bold R^3\Gamma^{(\eta_0)}_{\text{\rm Iw},h}(D,V)_{\Gamma}=0.
$
 This completes the proof of Theorem 4.2.3.

\proclaim{Corollary 4.2.5} The exponential map induces an isomorphism
of $D^{\Ph=p^{-1}}$ onto $\text{\rm coker} (f)_\Gamma \simeq \bold R^2\Gamma^{(\eta_0)}_{\text{\rm Iw},h}(D,V)_{\Gamma}$
and the diagram 
$$
\xymatrix{
D^{\Ph=p^{-1}} \ar[rr]^{\sim}   \ar[d]^{\lambda_D}& &\bold R^2\Gamma^{(\eta_0)}_{\text{\rm Iw},h}(D,V)^{\Gamma}
\ar[d]\\
D^{\Ph=p^{-1}} \ar[rr]^{(h-1)!\exp_V} & &\bold R^2\Gamma^{(\eta_0)}_{\text{\rm Iw},h}(D,V)_{\Gamma}}
$$
in which the map $\lambda_D$ is defined in Proposition 3.2.4, commutes.
\endproclaim

$\,$
\flushpar
{\bf 4.3. The module of $p$-adic $L$-functions.}
\newline
{\bf 4.3.1.} We conserve the notation and conventions of section 4.2. 
Let $D$ be an admissible subspace of $\Dc (V)$ and assume that $\Cal L(V,D)\ne 0.$
We review the definition of the module of $p$-adic $L$-functions using the formalism
of Selmer complexes.
Set 
$$
\Delta_{\text{\rm Iw},h} (D,V)\,=\,
{\det}^{-1}_{\Lambda_{\Bbb Q_p}}\left (\RG_{\text{\rm Iw},S}^{(\eta_0)}(V) \oplus \left (\underset{v\in S}\to \oplus
\RG_{\text{\rm Iw},f}^{(\eta_0)}(\Bbb Q_v,D,V)\right )\right ) \otimes 
{\det}_{\Lambda_{\Bbb Q_p}} \left (\underset{v\in S}\to \oplus
\RG_{\text{\rm Iw}}^{(\eta_0)}(\Bbb Q_v,V)\right ).
$$
The exact triangle
$$
\bold R\Gamma_{\text{\rm Iw},S}^{(\eta_0)}(D,V)@>>>
\left (\RG_{\text{\rm Iw},S}^{(\eta_0)}(V) \oplus \left (\underset{v\in S}\to \oplus
\RG_{\text{\rm Iw},f}^{(\eta_0)}(\Bbb Q_v,D,V)\right )\right )\otimes{\Cal H} @>>>
\left (\underset{v\in S}\to \oplus
\RG_{\text{\rm Iw}}^{(\eta_0)}(\Bbb Q_v,V)\right )\otimes \Cal H
$$
gives an isomorphism
$
\Delta_{\text{\rm Iw},h} (D,V)\otimes_{\Lambda_{\Bbb Q_p}}\Cal H \,\simeq\,
{\det}_{\Cal H}^{-1}\bold R\Gamma_{\text{\rm Iw},S}^{(\eta_0)}(D,V).
$
Let $\Cal K$ denote the  field of fractions of $\Cal H.$
By Theorem 4.2.3, all  $\bold R^i\Gamma_{\text{\rm Iw},S}^{(\eta_0)}(D,V)$
are $\Cal H$-torsion and we have a canonical map. 
$$
{\det}_{\Cal H}^{-1}\bold R\Gamma_{\text{\rm Iw},S}^{(\eta_0)}(D,V) \simeq
\underset{i\in\{2,3\}}\to \otimes {\det}^{(-1)^{i+1}}_{\Cal H}\bold R^i\Gamma_{\text{\rm Iw},S}^{(\eta_0)}(D,V)
\hookrightarrow \Cal K.
$$
The composition of these maps  gives a trivialization
$
i_{V,\text{\rm Iw},h}\,\,:\,\, 
\Delta_{\text{\rm Iw},h} (D,V) @>>> \Cal K.
$
 Fix a $\Bbb Z_p$-lattice $N$ of  $D$ and set
$$
\Delta_{\text{\rm Iw},h} (N,T)\,=\,
{\det}^{-1}_{\Lambda}\left (\RG_{\text{\rm Iw},S}^{(\eta_0)}(T) \oplus \left (\underset{v\in S}\to \oplus
\RG_{\text{\rm Iw},f}^{(\eta_0)}(\Bbb Q_v,N,T)\right )\right ) \otimes 
{\det}_{\Lambda} \left (\underset{v\in S}\to \oplus
\RG_{\text{\rm Iw}}^{(\eta_0)}(\Bbb Q_v,T)\right ).
$$

Perrin-Riou \cite{PR2} defined the module of $p$-adic $L$-functions associated to $(N,T)$  as
$$
\bold L_{\text{\rm Iw},h}^{(\eta_0)}(N,T)  =
i_{V,\text{\rm Iw},h} \left (\Delta_{\text{\rm Iw},h} (N,T)\right ) \subset \Cal K.
$$
Fix a generator $f(\gamma_1-1)$ of $\bold L_{\text{\rm Iw},h}^{(\eta_0)}(N,T)$ and define
a meromorphic $p$-adic function
$$
L_{\text{\rm Iw},h}(T,N,s)=f(\chi (\gamma)^s-1).
$$
\flushpar
Let $\omega_N$ be a generator  of ${\det}_{\Bbb Z_p}(N).$ The isomorphism $D\simeq t_V(\Bbb Q_p)$ allows us to consider 
$\omega_N$ as a basis of ${\det}_{\Bbb Q_p} t_V(\Bbb Q_p).$ We also fix a generator
$\omega_T$ of ${\det}_{\Bbb Z_p}T^+$ and define the $p$-adic period $\Omega_p (\omega_N,\omega_{\text{\rm T}})\in \Bbb Q_p$ 
by 
$
\omega_{\text{\rm B}}=\Omega_p (\omega_T,\omega_{\text{\rm B}})\omega_T.
$
Now we can state the main result of this paper.

\proclaim{Theorem 4.3.2} Assume that a pseudo-geometric representation  $V$ satisfies {\bf C1-5)}.  
Let $D$ be an admissible subspace of $\Dc (V).$ Fix a $G_{\Bbb Q}$-stable lattice $T$ of $V$ 
and a lattice $N$ of $D$. Assume that $\Cal L(D,V)\ne 0.$ Then

i)  $L_{\text{\rm Iw},h}(T,N,s)$ is a meromorphic $p$-adic function which has a zero 
at $s=0$ of order $e=\dim_{\Bbb Q_p}(D^{\Ph=p^{-1}}).$

ii) Let  $L_{\text{\rm Iw},h}^*(T,N,0)= \lim_{s\to 0} s^{-e}L_{\text{\rm Iw},h}(T,N,s)$
be the special value of $L_{\text{\rm Iw},h}(T,N,s)$ at $s=0.$ Then
$$
L_{\text{\rm Iw},h}^*(T,N,0)\overset{p}\to\sim \Gamma (h)^{d_{+}(V)} \,\Cal L (D,V)\, 
\,E_p^*(V,1)\,{\det}_{\Bbb Q_p} \left (\frac{1-p^{-1}\Ph^{-1}}{1-\Ph} \,\vert D_{-1} \right ) 
\frac{i_{\omega_N,\omega_{\text{\rm B}}, p}\,(\Delta_{\text{\rm EP}} (T))}{\Omega_p(\omega_T,\omega_{\text{\rm B}})},
$$
where $E_p(V,t)=E_p^*(V,t)\,\left (1-\dsize\frac{t}{p}\right )^e$ and $\Gamma (h)=(h-1)!.$
\endproclaim
\flushpar
{\bf 4.3.3. Proof of Theorem 4.3.2.} 
\newline
{\bf 4.3.3.1.} First recall  the formalism of Iwasawa descent which will be used in the proof.
The result we need is proved in \cite{BG}. This is a particular case of Nekov\'a\v r's
descent theory \cite{N2}. 
Let $C^\bullet$ be a perfect complex of $\Cal H$-modules and let $C^\bullet_0=C^\bullet\otimes^{\bold L}_{\Cal H}\Bbb Q_p.$
We have a natural distinguished triangle 
$$
C^\bullet @>X>>C^\bullet @>>>C^\bullet_0,
$$
where $X=\gamma_1-1.$ In each degree this triangle gives a short exact sequence
$$
0@>>>H^n(C^\bullet)_{\Gamma_1}@>>>H^n(C^\bullet_0)@>>>H^{n+1}(C^\bullet)^{\Gamma_1}@>>>0.
$$
One says that $C^\bullet$ is semisimple if the natural map
$$
H^n(C^\bullet)^{\Gamma_1}@>>>H^n(C^\bullet)@>>>H^n(C^\bullet)_{\Gamma_1}  \tag{4.7}
$$
is an isomorphism in all degrees. If $C^\bullet$ is semisimple, there exists a
natural trivialisation of ${\det}_{\Bbb Q_p}C^\bullet_0$, namely 
$$
\multline
\vartheta\,\,:\,\,{\det}_{\Bbb Q_p}C^\bullet_0 \simeq \underset{n\in \Bbb Z}\to\otimes
{\det}_{\Bbb Q_p}^{(-1)^n}H^n(C_0) \simeq 
\underset{n\in \Bbb Z}\to\otimes \left ({\det}^{(-1)^n}_{\Bbb Q_p}H^n(C^\bullet)_{\Gamma_1} \otimes 
{\det}^{(-1)^n}_{\Bbb Q_p}H^{n+1}(C^\bullet)^{\Gamma_1} \right ) \\
\simeq\underset{n\in \Bbb Z}\to\otimes \left ({\det}^{(-1)^n}_{\Bbb Q_p}H^n(C^\bullet)_{\Gamma_1} \otimes 
{\det}^{(-1)^{n-1}}_{\Bbb Q_p}H^{n}(C^\bullet)^{\Gamma_1} \right )
\simeq \Bbb Q_p
\endmultline
$$
where the last map is induced by (4.7). 
We now suppose that $C\otimes_{\Cal H}\Cal K$ is acyclic and write 
$
i_\infty \,\,:\,\, {\det}_{\Cal H} C^\bullet @>>> \Cal K
$
for the associated morphism in $\Cal P (\Cal K).$ Then
$
i_\infty ({\det}_{\Cal H} C^\bullet)=f\Cal H,
$
where $f\in \Cal K.$  Let $r$ be the unique integer such that 
$X^{-r}f$ is a unit of the localization $\Cal H_0$ of $\Cal H$ 
with respect to the principal ideal $X\Cal H$. 

\proclaim{Lemma 4.3.3.2} Assume that $C^\bullet$ is semisimple. Then 
$r=\dsize \underset{n\in \Bbb Z}\to \sum (-1)^{n+1} \dim_{\Bbb Q_p}H^n(C^\bullet)^{\Gamma_1}$
and there exists a commutative diagram
$$
\CD
{\det}_{\Cal H} C^\bullet @>X^{-r}i_\infty >>\Cal H_0\\
@V\otimes_{\Bbb Q_p}^{\bold L}VV  @VVV\\
{\det}_{\Bbb Q_p}{C^\bullet_0} @>\vartheta>> \Bbb Q_p
\endCD
$$
in which the right vertical arrow is the augmentation map.
\endproclaim
\demo{Proof} See \cite{BG}, Lemma 8.1. Remark that Burns and Greither consider  complexes
over $\Lambda \otimes_{\Bbb Z_p}  \Bbb Q_p$ but since $\Cal H$ is a B\'ezout ring,  all their
arguments work in our case and are omitted here. 
\enddemo
\flushpar
{\bf 4.3.3.3.}  By Theorem 4.2.3 the complex $\RG_{\text{Iw},h}^{(\eta_0)}(D,V)$
is semisimple and  the first assertion follows from Lemma 4.3.3.2 together with Corollary  4.2.5.
\flushpar
{\bf 4.3.3.4.}
In this subsection we compare the Bloch-Kato local condition at $p$ with the local
condition coming from Perrin-Riou's theory. Set $\RG_{f} (\Bbb Q_p,D,V)=D[-1]$ and define
$$
S=\text{\rm cone} \left (
\frac{1-p^{-1}\Ph^{-1}}{1-\Ph}\,:\, \RG_{f} (\Bbb Q_v,D,V)@>>>\RG_{f} (\Bbb Q_p,V)\right )\,[-1].
\tag{4.8}
$$
Thus, explicitly
$$
S=\left [ D\oplus \Dc (V)@>>>\Dc (V)\oplus t_V(\Bbb Q_p) \right ]\,[-1]\simeq 
\left [ D\oplus \Dc (V)@>>>\Dc (V)\oplus D) \right ]\,[-1],
$$
where the unique non-trivial map is given by
$$
(x,y)\mapsto \left ( (1-\Ph)\,y,\, \left (\frac{1-p^{-1}\Ph^{-1}}{1-\Ph}\,x+y\right )\,\pmod{\F^0\Dc (V)}
\right ).
$$
Thus
$
H^1(S)\,=\,D^{\Ph=p^{-1}}$ and $H^2(S)\,=\,\dsize\frac{t_V(\Bbb Q_p)}{(1-p^{-1}\Ph^{-1})D}\simeq
\dsize\frac{D}{(1-p^{-1}\Ph^{-1})D}$.
From the semi-simplicity of $\dsize\frac{1-p^{-1}\Ph^{-1}}{1-\Ph}$ it follows
that the natural projection  $H^1( S)@>>>H^2( S)$
is an  isomorphism and we have a  canonical trivialization
$
\vartheta_S\,:\,{\det}_{\Bbb Q_p}  S\simeq {\det}_{\Bbb Q_p}^{-1} H^1( S) \otimes 
{\det}_{\Bbb Q_p} H^2(S)\simeq \Bbb Q_p.
$
Hence the distingushed triangle 
$$
 S@>>> \RG_{f} (\Bbb Q_p,D,V) @>>>\RG_{f} (\Bbb Q_p,V)@>>>  S[1]
$$
 induces isomorphisms
$$
{\det}_{\Bbb Q_p} \RG_{f} (\Bbb Q_p,V) \simeq  \RG_{f} (\Bbb Q_p,D,V) \otimes
{\det}_{\Bbb Q_p}^{-1}   S \overset{\text{\rm via $\vartheta_S$}}\to\simeq  {\det}_{\Bbb Q_p} \RG_{f} (\Bbb Q_p,D,V).
$$

\proclaim{Lemma 4.3.3.5} i) Let $f\,:\,W@>>>W$ be a semi-simple endomorphism
of a finitely dimensional $k$-vector space $W$. The canonical
projection $\ker (f)@>>>\text{\rm coker} (f)$ is an isomorphism and 
the tautological
exact sequence 
$$
0@>>>\ker (f)@>>>W@>f>>W@>>>\text{\rm coker} (f)@>>>0
$$
induces an isomorphism
$$
{\det}^*f\,:\, {\det}_k (W) @>>> {\det}_k (W)\otimes  {\det}_k (\ker (f))\otimes
 {\det}^{-1}_k (\text{\rm coker} (f)) @>>>{\det}_k (W).
$$
Then
$
{\det}^*f (x)={\det} (f\,\vert\, \text{\rm coker} (f)).
$

ii) The diagram
$$
\xymatrix{
{\det}_{\Bbb Q_p} \RG_{f} (\Bbb Q_p,V) \ar[rr] \ar[d] & &{\det}_{\Bbb Q_p} \RG_{f} (\Bbb Q_p,D,V) \ar [dd]\\
{\det}_{\Bbb Q_p}^{-1} t_V(\Bbb Q_p) \ar[d] & &\\
{\det}_{\Bbb Q_p}^{-1} D \ar[rr]^{\det^*\left (\frac{1-p^{-1}\Ph^{-1}}{1-\Ph}\vert D\right )E_p(V,1)} & & {\det}_{\Bbb Q_p}^{-1}D}
$$
in which the bottom map is the multiplication by $\det^*\left (\frac{1-p^{-1}\Ph^{-1}}{1-\Ph} \vert D \right ) 
  E_p(V,1)$,
commutes.
\endproclaim
\demo{Proof} The proof of i) is straightforward and is omitted here. Next, ii) follows from i)
applyed to $W=D$ and the fact what $E_p(V,1)=\det \left (1-\Ph\,\vert \,\Dc (V) \right ).$
\enddemo

{\bf 4.3.3.6.} Now we can prove Theorem 4.3.2. Define

$$
\RG_{f} (\Bbb Q_v,N,T)= \RG_{\text{\rm Iw},f}^{(\eta_0)}(\Bbb Q_v,N,T)\otimes_{\Lambda}^{\bold L} \Bbb Z_p,
\qquad  \RG_{f} (\Bbb Q_v,D,V)=\RG_{f} (\Bbb Q_v,N,T)\otimes_{\Bbb Z_p}\Bbb Q_p.
$$
Remark that for $v=p$ this definition coincides with the definition given in 4.3.3.4.
Applying $\otimes^{\bold L}_{\Cal H}\Bbb Q_p$ to   the map  
$\RG_{\text{\rm Iw},f}^{(\eta_0)} (\Bbb Q_v,D,V) @>>> \RG_{\text{\rm Iw}}^{(\eta_0)}(\Bbb Q_v,T)
\otimes^{\bold L}_{\Lambda} \Cal H$ we obtain a morphism
$$
\RG_{f} (\Bbb Q_v,D,V)@>>> \RG (\Bbb Q_v, V).
$$ 
If $v\ne p,$ then $\RG_{f} (\Bbb Q_v,D,V)=\RG_f(\Bbb Q_v,V)$ and this morphism coincides
with the natural map $\RG_f(\Bbb Q_v,V)@>>> \RG (\Bbb Q_v,V).$
If $v=p,$ then $\RG_{f} (\Bbb Q_v,D,V)=D[-1]$ and by Theorem 2.2.4 it coincides with
the composition
$$
D@>\frac{1-p^{-1}\Ph^{-1}}{1-\Ph}>> \Dc (V)@>(h-1)! \exp_{V,\Bbb Q_p}>> H^1(\Bbb Q_p,V).
$$

Let $\RG_{f,h} (D,V)$ denote the Selmer complex associated to the diagram
$$
\xymatrix{
\RG_{S} (V)
\ar[r] & \underset{v\in S} \to  \oplus \RG (\Bbb Q_v,V)\\
 & \underset{v\in S}\to \oplus \RG_{f}(\Bbb Q_v,D,V)
 \ar[u]
}
$$

Then we have a distinguished triangle
$$
\RG_{f,h} (D,V) @>>>\RG_{S} (V)\oplus \left (\underset{v\in S}\to \oplus \RG_{f}(\Bbb Q_v,D,V)\right )@>>>\underset{v\in S} \to  \oplus \RG (\Bbb Q_v,V)
\tag{4.9}
$$
which induces isomorphisms
$$
\align
&{\det}^{-1}_{\Bbb Q_p}\RG_S(V) \otimes_{\Bbb Q_p}\left (\underset{v\in S}\to \otimes {\det}_{\Bbb Q_p}\RG (\Bbb Q_v,V)\right )
\otimes {\det}_{\Bbb Q_p}D \iso {\det}^{-1}_{\Bbb Q_p}\RG_{f,h} (D,V),\\
& \xi_{D,h} \,:\, \Delta_{\text{\rm EP}}(V)\otimes_{\Bbb Q_p} \left (  {\det}_{\Bbb Q_p}D\otimes {\det}^{-1}_{\Bbb Q_p}V^+\right )\iso 
{\det}^{-1}_{\Bbb Q_p}\RG_{f,h} (D,V).
\endalign
$$

Next, $\RG_{f,h} (D,V)=\RG_{\text{\rm Iw},h}^{(\eta_0)} (D,V) \otimes_{\Cal H}\Bbb Q_p$
and for any $i$ one has an exact sequence
$$
0@>>>\bold R^i\Gamma_{\text{\rm Iw},h}^{(\eta_0)} (D,V)_{\Gamma} @>>>\bold R^i\Gamma_{f,h} (D,V)
@>>>\bold R^{i+1}\Gamma_{\text{\rm Iw},h}^{(\eta_0)} (D,V)^{\Gamma}@>>>0.
$$
From Theorem 4.2.3 it follows that 
$$
\bold R^i\Gamma_{f,h} (D,V)=\cases \bold R^{2}\Gamma_{\text{\rm Iw},h}^{(\eta_0)} (D,V)^{\Gamma}
&{\text{\rm if $i=1$}}\\
\bold R^{2}\Gamma_{\text{\rm Iw},h}^{(\eta_0)} (D,V)_{\Gamma} &{\text{\rm if $i=2$}}\\
 0 &{\text{\rm if $i\ne 1,2$}}.
\endcases
$$
Therefore, the   isomorphism $\bold R^2\Gamma_{\text{\rm Iw},h}(D,V)^{\Gamma}@>>>
\bold R^2\Gamma_{\text{\rm Iw},h}(D,V)_{\Gamma}$ induces a canonical trivialization
$$
\vartheta_{D,h}\,:\,{\det}_{\Bbb Q_p}\RG_{f,h} (D,V)  \iso  \Bbb Q_p.
$$
 By Lemma 4.3.3.2
we have a commutative diagram
$$
\CD
{\det}_{\Cal H}^{-1}\RG_{\text{\rm Iw},h}^{(\eta_0)} (D,V) @>X^{-e}i_{V,\text{\rm Iw},h}>>
{\Cal H}_0
\\
@V\,^{\bold L}\otimes_{\Cal H}\Bbb Q_p VV    @VVV \\
{\det}_{\Bbb Q_p}^{-1}\RG_{f,h} (D,V)  @>\vartheta_{D,h}^{-1}>>  \Bbb Q_p.
\endCD
$$
Since 
$$
\Delta_{\text{\rm Iw},h}(N,T)\otimes^{\bold L}_{\Lambda}\Bbb Z_p \simeq \Delta_{\text{\rm EP}}(T)\otimes_{\Bbb Z_p}
\omega_N \otimes_{\Bbb Z_p} \omega^{-1}_T 
$$ 
it implies that 
$$
\vartheta_{D,h}^{-1}\circ  \xi_{D,h} (\Delta_{\text{\rm EP}}(T)\otimes_{\Bbb Z_p}
\omega_N \otimes_{\Bbb Z_p} \omega^{-1}_T )    \,=\,
\log (\chi (\gamma))^{-e} L^*_{\text{\rm Iw},h}(T,N,0)\,\Bbb Z_p.
\tag{4.10}
$$

Consider the  diagram
$$
\xymatrix{
\RG_f(V) \ar[r] &\RG_S(V)\oplus \underset{v\in S\cup\{\infty\}}\to \oplus \RG_f(\Bbb Q_v,V) \ar[r]&
 \underset{v\in S\cup\{\infty\}}\to \oplus \RG(\Bbb Q_v,V)\\
\RG_{f,h}(D,V) \ar[r]\ar[u] &\RG_S(V)\oplus \underset{v\in S}\to \oplus \RG_f(\Bbb Q_v,D,V) \ar[r] \ar[u]
&\underset{v\in S}\to \oplus \RG(\Bbb Q_v,V) \ar[u]\\
L \ar[r]\ar[u] &S\oplus V^+[-1] \ar[r] \ar[u] &V^+[-1] \ar[u]
}
\tag{4.11}
$$ 
in which  $L=\text{\rm cone} \left (\RG_{f,h}(D,V) @>>> \RG_f(V) \right )\,[-1]$ and the upper
and  middle rows coincide with (4.1) and (4.9) up to the following modification:
the  map ${\text{\rm loc}}_p\,:\, \RG_f(\Bbb Q_p,V)@>>> \RG (\Bbb Q_p,V)$ 
is replaced by
$\Gamma (h)\,{\text{\rm loc}}_p$.
It follows from {\bf C1-5)}  that $\RG_f(V)$ is acyclic. Hence in the derived
category $\Cal D^p(\Bbb Q_p)$ the composition
$
\alpha \,:\,S\iso L\iso \RG_{f,h}(D,V)
$
is an isomorphism.
An easy diagram search shows that $H^1(S)\simeq \bold R^1 \Gamma_{f,h}(D,V)$ coincides
with ${\text{\rm id}}\,:\,D^{\Ph=p^{-1}} @>>> D^{\Ph=p^{-1}}$ and that 
 $H^2(S)\simeq \bold R^2 \Gamma_{f,h}(D,V)$ coincides with $\Gamma (h)\exp_{V,\Bbb Q_p} .$
Therefore, we have a commutative diagram 
$$
\CD
{\det}_{\Bbb Q_p} S  @> \alpha>>  {\det}_{\Bbb Q_p}\RG_{f,h} (D,V)\\
@V{\vartheta_S}VV  @V{\vartheta_{D,h}}VV\\
\Bbb Q_p @>\kappa>> \Bbb Q_p
\endCD
$$

there $\kappa$  can be written as the composition
$$
\Bbb Q_p \iso {\det}^{-1}_{\Bbb Q_p} H^1(S)\otimes   {\det}_{\Bbb Q_p} H^2(S)
\iso   {\det}^{-1}_{\Bbb Q_p}\bold R^1 \Gamma_{f,h}(D,V) \otimes
 {\det}_{\Bbb Q_p}\bold R^2 \Gamma_{f,h}(D,V) \iso \Bbb Q_p
$$
From Proposition 3.2.4  and  Corollary 4.2.5 we obtain immediately that 
$$
\kappa\,=\,(\log \chi (\gamma))^{e} \left (1-\dsize\frac{1}{p}\right )^{e}
\Cal L(D,V)^{-1}\, {\text{\rm id}}_{\Bbb Q_p}.
\tag{4.12}
$$
Passing to determinants in the diagram (4.11) we obtain  a commutative diagram
$$
\xymatrix{
\Delta_{\text{\rm EP}} (V) \otimes 
\left ({\det}_{} (t_V(\Bbb Q_p))\otimes {\det}^{-1}_{}V^+ \right )\ar[d]^{f} \ar[r]   &{\det}_{}^{-1}
\RG_f (V) \ar[r]  \ar[d] & \Bbb Q_p \ar @{=}[d]\\
\Delta_{\text{\rm EP}} (V) \otimes 
\left ({\det}_{} D \otimes {\det}^{-1}_{}V^+ \right ) \otimes {\det}_{}S \ar[r]^{\xi_{D,h}\otimes 
\alpha}  \ar[d]^{\text{\rm id}\otimes \vartheta_S}&{\det}_{}^{-1}
\RG_{f,h} (D,V) \otimes {\det}_{}
\RG_{f,h} (D,V)  \ar[r]^-{\text{\rm duality}} \ar[d]^{{\text{\rm id}}\otimes \vartheta_{D,h}}& \Bbb Q_p \ar @{=}[d]\\
\Delta_{\text{\rm EP}}^{} (V) \otimes 
\left ({\det}_{} D \otimes {\det}^{-1}_{\Bbb Q_p}V^+ \right )
\ar[r]^-{\xi_{D,h}\otimes \kappa} &{\det}_{}^{-1}
\RG_{f,h} (D,V) \ar[r]^-{\vartheta_{D,h}^{-1}}    &\Bbb Q_p}
$$
in which the map $f$ is induced by (4.8). 
The upper row of this diagram sends $\Delta_{EP}(T)\otimes  (\omega_N\otimes \omega^{-1}_{\text{\rm B}})$ onto 
$$\Gamma (h)^{d_{+}(V)} i_{\omega_N,\omega_{\text{\rm B}},p} (\Delta_{EP}(T)).
\tag{4.13}$$
 From Lemma 4.3.3.5  it follows that 
 $$(\text{\rm id}\otimes \vartheta_S)\circ f\,=\,{\det}^* \left (\frac{1-p^{-1}\Ph^{-1}}{1-\Ph}\mid D \right )^{-1} E_p(V,1)^{-1}\,\text{\rm id}
 \tag{4.14}
$$
Next,  (4.10) and (4.12) give 
$$
\vartheta_{D,h}^{-1}\circ (\xi_{D,h}\otimes \kappa) (\Delta_{\text{\rm EP}}(T)\otimes_{\Bbb Z_p} \omega_N \otimes_{\Bbb Z_p} \omega^{-1}_T)\,=\, 
\left (1-\dsize\frac{1}{p}\right )^{e}
\Cal L(D,V)^{-1}      L^*_{\text{\rm Iw},h}(T,N,0)\,\Bbb Z_p.
\tag{4.15}
$$

Putting together (4.13), (4.14) and (4.15) we obtain that
$$
L_{\text{\rm Iw},h}^*(T,N,0)\overset{p}\to\sim \Gamma (h)^{d_{+}(V)} \,\Cal L (D,V)\, 
\,E_p^*(V,1)\,{\det}^*_{\Bbb Q_p} \left (\frac{1-p^{-1}\Ph^{-1}}{1-\Ph} \,\vert D \right ) 
\frac{i_{\omega_N,\omega_{\text{\rm B}}, p}\,(\Delta_{\text{\rm EP}} (T))}{\Omega_p(\omega_T,\omega_{\text{\rm B}})}.
$$
The theorem is proved.

\head  {\bf Appendix. Galois cohomology of $p$-adic
representations}
\endhead

{\bf A.1.} Let $K$ be a finite extension of $\Bbb Q_p$ and $T$  a
$p$-adic representation of $G_K.$ Fix a topological generator $\gamma$ of $\Gamma$.
Let $\bD (T)=(T\otimes_{\Bbb Z_p} \bold A)^{H_K}$ be the
$(\Ph,\Gamma)$-module associated to $T$ by Fontaine's theory \cite{F2}. Consider the complex
$$
C_{\Ph,\g}(\bD (T))\,=\,\left [\bD(T)@>f>> \bD (T)\oplus \bD (T) @>g>>
\bD (T) \right ]
$$
where the  modules are placed in degrees $0$, $1$ and $2$ and the
maps $f$ and $g$ are given by
$$
f(x)= ((\Ph-1)\,x\,,\, (\g-1)\,x),\qquad g(y,z)= (\g-1)\,y\,-\,
(\Ph-1)\,z.
$$

\proclaim{Proposition A.2} There are canonical and functorial
isomorphisms
$$
h^i\,:\, H^i(C_{\Ph,\g}(\bD (T))) \iso H^i(K,T)
$$
which can be described explicitly by the following formulas:

i) If $i=0,$ then $h^0$ coincides with the natural isomorphism
$$
\bD (T)^{\Ph=1,\g=1}\,=\,H^0(K,T\otimes_{\Bbb Z_p}\A^{\Ph=1})\,=\,
H^0(K,T).
$$
ii)  Let $\alpha,\beta \in \bD (T)$ be such that $(\g-1)\,\alpha
\,=\,(1-\Ph)\,\beta.$ Then $h^1$ sends $\cl (\alpha,\beta)$ to the
class of the cocycle
$$
\mu_1 (g)\,=\, (g-1)\,x\,+\,\dsize \frac{g-1}{\g-1}\,\beta,
$$
where $x \in \bD (T)\otimes_{\A_K} \A$ is a solution of the
equation $(1-\Ph)\,x\,=\,\alpha .$

iii) Let $\widehat\g\in G_K$ be a lifting of $g\in \Gam$ and let
$x$ be a solution of $(\Ph -1)\,x\,=\,\alpha .$ Then  $h^2$ sends
$\alpha$ to the class of the 2-cocycle
$$
\mu_2 (g_1,g_2)\,=\,\widehat \g^{k_1}(h_1-1)\,\frac{\widehat
\g^{k_2}-1}{\widehat \g-1}\,x
$$
where $g_i\,=\,\hat\g^{k_i}h_i,$ $h_i\in H_K.$
\endproclaim
\demo{Proof} The isomorphisms $h^i$ were constructed   in
\cite{H1}, Theorem 2.1. Remark that i) follows directly from this
construction (see \cite{H1}, p.573) and that ii) is proved in
\cite{Ben1}, Proposition 1.3.2 and \cite{CC2}, Proposition I.4.1.
The proof of iii) follows along exactly the same lines. Namely, it
is enough to prove this formula modulo $p^n$ for each $n$. Let $\alpha \in \bD
(T)/p^n\bD (T) .$ By Proposition 2.4 of \cite{H1} there exists
$r\ge 0$ and $y\in \bD(T)/p^n\bD (T)$ such that
$(\Ph-1)\,\alpha\,=\,(\g-1)^r\beta .$ Let
$$
N_x=(\bD (T)/p^n\bD (T))\oplus (\oplus_{i=1}^r
(\A_K/p^n\A_K)\,t_i),
$$
where $\Ph (t_i)\,=\,t_i+(\g-1)^{r-i}(\alpha)$ and $\g
(t_i)\,=\,t_i+t_{i-1}.$ Then $N_x$ is a $(\Ph,\Gamma)$-module and
we have a short exact sequence
$$
0@>>>\bD@>>>N_x@>>>X@>>>0
$$
where $X=N_x/M\simeq \oplus_{i=1}^r \A_K/p^n\A_K \bar t_i.$ An
easy diagram search shows that the connecting homomorphism
$\delta^1_{\bD}\,:\,H^1(C_{\Ph,\g}(\bD (X))) @>>>H^2(C_{\Ph,\g}(\bD(T)))$
sends $\cl (0,\bar t_r)$ to $-\cl (\alpha).$ The functor 
$\bold V(D)=(D\otimes_{\bold A_K} \bold A)^{\Ph=1}$ is a quasi-inverse 
to $\bD $. Thus one has an exact sequence of Galois
modules
$$
0@>>> T/p^nT@>>> T_x @>>> \bold V(X)@>>>0
$$
where $T_x=\bold V(N_x).$ From the definition of $x$ it follows
immediately that $t_r-x\in T_x .$ By ii), $h^1(\cl (0,\bar t_r))$
can be represented by the cocycle $c(g)\,=\,
\dsize\frac{g-1}{\g-1}\,\bar t_r$ and we fix its lifting $\hat
c\,:\,G_K@>>>N_x$ putting $\hat c(g)\,=\, \dsize\frac{g-1}{\g-1}\,
(t_r-x).$ As $g_1\hat c(g_2)-\hat c(g_1g_2)+\hat c(g_1)\,=\, -\mu_2(g_1,g_2), $
the connecting map $\delta_T^1\,:\,H^1(K,\bold
V(X))@>>>H^2(K,T/p^nT)$ sends $\cl (c)$ to $-\cl (\mu_2)$ and iii)
follows from    the commutativity of the diagram
$$
\CD
H^1(C_{\Ph,\g}(X)) @>\delta^1_{\bD}>>H^2(C_{\Ph,\g}(T/p^nT))\\
@Vh^1VV @Vh^2VV\\
H^1(K,\bold V(X))@>\delta^1_T>>H^2(K,T/p^nT).
\endCD
$$
\enddemo




\proclaim{Proposition A.3} The complexes $\RG (K,T)$ and  $C_{\Ph,\g}(T)$ 
are isomorphic in  $\Cal D(\Bbb Z_p)$. 


\endproclaim
\demo{Proof} 
The proof  is standard (see for example \cite{BF},
proof of Proposition 1.17). The exact
sequence
$$
0@>>> T@>>>\bD (T)\otimes_{\A_K}\A@>\Ph-1>>\bD
(T)\otimes_{\A_K}\A@>>>0
$$
gives rise to an exact sequence of complexes
$$
0@>>> C^{\bullet}_{\co}(G_K,T)@>>>C^{\bullet}_{\co}(G_K,\bD
(T)\otimes_{\A_K}\A)@>\Ph-1>>C^{\bullet}_{\co} (G_K,\bD
(T)\otimes_{\A_K}\A)@>>>0
$$
Thus $\bold R\Gamma (K,T)$ is quasi-isomorphic to the total
complex
$$
K^\bullet (T)\,=\,\text{Tot}^{\bullet} \left
(C^{\bullet}_{\co}(G_K,\bD
(T)\otimes_{\A_K}\A)@>\Ph-1>>C^{\bullet}_{\co} (G_K,\bD
(T)\otimes_{\A_K}\A)\right ).
$$
On the other hand
$
C_{\Ph,\g}(T)\,=\,\text{Tot}^{\bullet} \left
(A^{\bullet}(T)@>\Ph-1
>>A^{\bullet}(T) \right ),
$
where $A^{\bullet}(T)=[\bD (T) @>\gamma-1>>\bD (T)]$.
Consider the following commutative diagram of complexes
$$
\xymatrix{ \bD (T) \ar[d]^{\beta_0} \ar[r]^{\g-1} &\bD (T) \ar[d]^{\beta_1} \ar[r] &0 \ar[r]\ar[d] & \cdots \\
C^{0}(G_K,\bD (T)\otimes_{\A_K}\A) \ar[r]& C^{1}(G_K,\bD
(T)\otimes_{\A_K}\A) \ar[r]& C^{2}(G_K,\bD (T)\otimes_{\A_K}\A)
\ar[r]& \cdots }
$$
in which  $\beta_0 (x)=x$ viewed as a constant function on $G_K$ and
$\beta_1(x)$ denotes the map $G_K@>>>\bD (T)\otimes_{\A_K}\A)$
defined by
$
(\beta_1 (x))\,(g)=\dsize\frac{g-1}{\g-1}\, x.
$
This diagram induces a map $\Tot^{\bullet}
(A^\bullet(T)@>\Ph-1>>A^{\bullet}(T)) @>>>K^\bullet (T)$ and we
obtain a diagram
$$
C_{\Ph,\g}(T) @>>>K^{\bullet}(T) \leftarrow \RG (K,T)
$$
where the right map is a quasi-isomorphism. Then for each $i$ one
has a map
$$
H^i(C_{\Ph,\g}(T)) @>>>H^i(K^{\bullet}(T)) \simeq H^i(K,T)
$$
and an easy diagram search shows that it coincides with $h^i.$ The
proposition is proved.
\enddemo
\proclaim{Corollary A.4} Let $V$ be a $p$-adic representation of $G_K$.
Then the complexes $\bold R\Gamma (K,V)$, $C_{\Ph,\g}(\bD^{\dag}(V))$ and
$C_{\Ph,\g}(\Ddagrig (V))$ are isomorphic in $\Cal D(\Bbb Q_p).$
\endproclaim
\demo{Proof} This follows from Theorem 1.1 of \cite{Li} together with Proposition A.2.
\enddemo 

{\bf A.5.} Recall that $K_\infty/K$ denotes the cyclotomic
extension obtained by adjoining all $p^n$-th roots of unity. Let
$\Gamma =\G (K_\infty/K)$ and let $\Lambda (\Gamma)\,=\,\Bbb Z_p[[\Gamma]]$
denote the Iwasawa algebra of $\Gam .$ For any  $\Bbb Z_p$-adic
representation  $T$ of $G_K$ the induced representation
$\Ind_{K_{\infty}/K}T$ is isomorphic to $(T\otimes_{\Bbb Z_p}\Lambda
(\Gamma))^\iota $ and we set
$
\R\Gamma_{\text {\rm Iw}} (K,T)\,=\,C^{\bullet}_{\co}
(G_K\,,\Ind_{K_{\infty}/K}T).
$
Consider the complex
$$
C_{\Iw, \psi}(T)\,=\,\left [\bD (T)@>\psi-1>>\bD (T)\right ]
$$
in which the first term is placed in degree $1$.

\proclaim{Proposition A.6} There are canonical and functorial
isomorphisms
$$
h^i_{\Iw}\,\,:\,\,H^i(C_{\Iw,\psi}(T))@>>> H^i_{\Iw}(K,T)
$$
which can be described explicitly by the following formulas:

i) Let $\alpha \in \bD (T)^{\psi=1}.$ Then $(\Ph -1)\,\alpha \in
\bD (T)^{\psi=0}$ and for any $n$ there exists a unique
$\beta_n\in \bD (T)$ such that $(\gn -1)\,\beta_n= (\Ph-1)\,\alpha
.$ The map $h^1_{\Iw}$ sends $\cl (\alpha)$ to
$
(h^1_n (\cl (\beta_n,\alpha)))_{n\in \bold N} \in H^1_{\Iw}
(K_n,T).
$

ii) If $\alpha \in \bD (T),$ then
$
h^2_{\Iw} (\cl (\alpha))\,=\,-(h^2_n(\Ph (\alpha)))_{n\in \bold
N}.
$
\endproclaim
\demo{Proof} The proposition follows from Theorem II.1.3 and
Remark II.3.2 of \cite{CC2} together with Proposition A.2.
\enddemo

\proclaim{Proposition A.7} The complexes $\R\Gamma_{\text {\rm
Iw}} (K,T)$ and $C_{\Iw, \psi}(T)$ are isomorphic in the derived
category $\Cal D (\Lambda (\Gamma)).$
\endproclaim
\demo{Proof} We repeat the arguments used in the proof of
Proposition A.1.2 with some modifications. For any $n\geqslant 1$ one
has an exact sequence
$$
0@>>>\Ind_{K_{n}/K}T @>>> (\bD(T)\otimes_{\Bbb Z_p}\Bbb
Z_p[G_n]^{\iota})\otimes_{\A_K}\A @>\Ph-1>> (\bD(T)\otimes_{\Bbb
Z_p}\Bbb Z_p[G_n]^{\iota})\otimes_{\A_K}\A @>>> 0.
$$
Set $\bD (\Ind_{K_\infty/K}T)=\bD (T)\otimes_{\Bbb Z_p}\Lambda
(\Gamma)^{\iota}$ and
$$
\bD (\Ind_{K_\infty/K}(T)) \hat \otimes_{\A_K}
\A\,=\,\varprojlim_n (\bD (T)\otimes_{\Bbb Z_p}\Bbb
Z_p[G_n]^{\iota})\otimes_{\A_K} \A.
$$
As $\Ind_{K_{n}/K}T$ are compact, taking projective limit one
obtains an exact sequence
$$
0@>>>\Ind_{K_{\infty}/K}T @>>> \bD (\Ind_{K_\infty/K}(T)) \hat
\otimes_{\A_K} \A @>\Ph-1>> \bD (\Ind_{K_\infty/K}(T)) \hat
\otimes_{\A_K} \A @>>>0 .
$$

Thus $\R_{\text {\rm Iw}} (K,T)$ is quasi-isomorphic to
$$
K^{\bullet}_{\Iw}(T)\,=\, \text{Tot}^{\bullet} \left
(C^{\bullet}_{\co}(G_K,\bD (\Ind_{K_\infty/K}T) \hat
\otimes_{\A_K} \A) @>\Ph-1>> C^{\bullet}_{\co}(G_K,\bD
(\Ind_{K_\infty/K}T) \hat \otimes_{\A_K} \A) \right ).
$$
We construct a quasi-isomorphism
$f_{\bullet}\,:\,C_{\Iw,\psi}(T)@>>> K^{\bullet}_{\Iw}(T).$ Any
$x\in \bD (T)$  can be written in the form $x=(1-\Ph
\psi)\,x\,+\,\Ph \psi (x)$ where $\psi (1-\Ph \psi)\,x\,=\,0.$
Then for each $n \geqslant  0$ the equation
$
(\gn -1)\,y_n\,=\,(\Ph \psi -1)\,x
$
has a unique solution  $y_n\in \bD (T)^{\psi=0}$ (\cite{CC2},
Proposition I.5.1). In particular,
$
y_n\,=\,\dsize \frac{\gamma_{n+1}-1}{\gn-1}\,y_{n+1}
$
and we have a compatible system of elements
\linebreak
$
Y_n=\dsize\sum_{k=0}^{|G_n|-1} \g^k \otimes \g^k(y_n)\,\in \bD (T)
\otimes_{\Bbb Z_p} \Bbb Z_p[G_n]^{\iota}.
$
Put $Y=(Y_n)_{n\ge 0}  \in \bD (\Ind_{K_\infty/K}T).$ Then
$$
(\gn-1)\,Y_n\,=\,(\g-1)\,Y \pmod{\bD (\Ind_{K_n/K}T)}.
$$
Let $\eta_x \in C^{1}_{\co}(G_K,\bD (\Ind_{K_\infty/K}T)\hat
\otimes_{\A_K} \A)$ be the map defined by
$
\eta_x (g)\,=\,\dsize \frac{g-1}{\g-1}\,(1\otimes x).
$
Define
$
f_1\,:\,\bD (T)@>>> K^1_{\Iw}(T)\,=\, C^{0}_{\co}(G_K,\bD
(\Ind_{K_\infty/K}T)\hat \otimes_{\A_K} \A) \oplus
C^{1}_{\co}(G_K,\bD (\Ind_{K_\infty/K}T)\hat \otimes_{\A_K} \A)
$
by $f_1(x)=(Y,\eta_x)$ and
$
f_2\,:\,\bD (T) @>>> C^{1}_{\co}(G_K,\bD (\Ind_{K_\infty/K}T)\hat
\otimes_{\A_K} \A)\subset K^2_{\Iw}(T)
$
by $f_2(z)= -\eta_{\Ph (z)}.$ 
It is easy to check that $f_{\bullet}$ is a morphism of complexes. This gives   a diagram
$$
C_{\Iw,\psi}(T) @>>>K^{\bullet}_{\Iw}(T) \leftarrow \RG_{\Iw}
(K,T)
$$
in which the right map is a quasi-isomorphism. Using Proposition A.1.4 it is not difficult to check 
that for each  $i$ the induced  map
$$
H^i(C_{\Iw,\psi}(T)) @>>>H^i(K^{\bullet}_{\Iw}(T)) \simeq
H_{\Iw}^i(K,T)
$$
coincides with $h^i_{\Iw}.$ The proposition is proved.
\enddemo

\proclaim{Corollary A.8} The complexes $\bold R\Gamma_{\Iw} (K,T)$ and
$C_{\Iw,\psi}^{\dag}(T)$ are isomorphic in $\Cal D(\Lambda (\Gamma)).$
\endproclaim
\demo{Proof} One has $\bD^{\dag} (T)^{\psi=1}=\bD (T)^{\psi=1}$ (\cite{CC1}, Proposition 3.3.2)
and $\bD^{\dag} (T)/(\psi-1) =\bD (T)/(\psi-1)$ (\cite{Li}, Lemma 3.6).
This shows that the inclusion $C_{\Iw,\psi}^{\dag}(T) @>>> \bD (T)^{\psi=1}$ is a quasi-isomorphism.
\enddemo

\flushpar 
{\bf Remark A.9.} These results can be slightly improved.
Namely, set $r_n=(p-1)p^{n-1}.$ The method used in the proof of Proposition III.2.1 
\cite{CC2} allows to show that $\psi (\bD^{\dag,r_n}(T))\subset \bD^{\dag,r_{n-1}}(T)$
for $n\gg 0.$ Moreover, for any $a\in  \bD^{\dag,r_n}(T)$ the solutions of 
the equation $(\psi-1)\,x=a$ are in $\bD^{\dag,r_n}(T).$ Thus
$C_{\Iw}^{\dagger,r_n} (T)\,=\, \left [\bD^{\dag,r_n}
(T)@>\psi-1>>\bD^{\dag, r_n} (T)\right ]$, $n\gg 0$ is a well-defined complex 
which is quasi-isomorphic to $C_{\Iw,\psi}^{\dag}(T).$  
Further, as $\Ph (\A^{\dagger,r/p})\,=\,\A^{\dagger,r}$ 
we can consider the complex 
$$
C^{\dagger ,r_n}_{\Ph,\g}(T)\,=\,\bigl [\bD^{\dagger
,r_{n-1}}(T)@>f>> \bD^{\dag,r_n} (T)\oplus \bD^{\dagger , r_{n-1}}
(T) @>g>> \bD^{\dag,r_n}
(T) \bigr ], \qquad n\gg 0
$$
in which $f$ and $g$ are defined by the same formulas as before.
Then the inclusion 
$
C^{\dagger, r_n}_{\Ph,\g} (T) @>>>C_{\Ph,\g}(T)
$
is a quasi-isomorphism.

\enddocument
\bye